\documentclass[a4paper,10pt]{article} \usepackage{amsmath,amssymb,amsthm}
\usepackage[english]{babel}

\newtheorem{theo}{Theorem}[section] \newtheorem{defi}[theo]{Definition}
\newtheorem{lemm}[theo]{Lemma} \newtheorem{prop}[theo]{Proposition}
\newtheorem{coro}[theo]{Corollary}

\newcommand{\Na}{\mathbb N}                   % N des entiers

\newcommand{\Ra}{\mathbb R}

 \newcommand{\R}{\mathbb R}                  % R des reels

\newcommand{\Ca}{\mathbb C}                   % C des complexes

\newcommand{\scal}[1]{\langle #1 \rangle}

\newcommand{\finpreuve}{\hfill $\Box$}

\newcommand{\name}{$\underline{\qquad \qquad}$}

\begin{document}

\author{  Jean-Marc
Bouclet\footnote{Jean-Marc.Bouclet@math.univ-lille1.fr}\\
 \\ Vincent Bruneau\footnote{Vincent.Bruneau@math.u-bordeaux1.fr} }
\title{{\bf \sc Semiclassical Resonances of Schr\"odinger operators as zeroes of regularized determinants}}

\maketitle

\begin{abstract} We prove that the resonances of long range perturbations of the (semiclassical) Laplacian  are the zeroes of  natural perturbation determinants. We more precisely obtain factorizations of these determinants of the form $ \prod_{w = {\rm resonances}}(z-w) \exp (\varphi_p(z,h)) $ and give semiclassical bounds on $ \partial_z \varphi_p $ as well as a representation of  Koplienko's regularized spectral shift function.  Here the index $ p \geq 1 $ depends on the decay rate at infinity of the perturbation.
\end{abstract}

\section{Introduction and results}
\setcounter{equation}{0}

One of the main purposes of Scattering Theory is the study of
selfadjoint operators with absolutely continuous (AC) spectrum. This
corresponds physically to extended or delocalized states, by
opposition to the localized or confined states which give rise to
discrete spectrum.
A typical mathematical example of confining system is given by
the Laplacian $ \Delta_g $  (or more general elliptic operators)
on a compact riemannian manifold: here, the states (ie the
eigenfunctions) are clearly localized by the compactness
assumption and the spectrum is a non decreasing sequence of
eigenvalues tending to infinity. 

Quite naively, $ \Delta_g $ can
 be viewed as an infinite dimensional analogue of an hermitian
matrix $ A = A^* $ on $ \Ca^N $. In that case, the spectrum of
$A$ is given by the roots of the characteristic polynomial $
\mbox{Det}(A-z) $.
% or, equivalently, by the poles of$ \mbox{det}(A-z)^{-1} $.
It is elementary to check that, for $ z
$ in the upper half plane,
\begin{eqnarray}
\mbox{Det}(A-z) = \exp \left( \partial_s \mbox{tr}
(A-z)^{s}_{|s=0} \right), \label{zetaelementaire}
\end{eqnarray}
so $ \mbox{Det}(A-z) $ can be defined as the analytic continuation
(with respect to $z$) of the right hand side of
(\ref{zetaelementaire})
 to the complex plane.
This is an elementary version of the classical definition of
determinants via a Zeta function (here $ \mbox{tr} (A-z)^{s} $),
which is used  in infinite dimension, typically for elliptic
operators on compact manifolds as initially introduced by Ray and Singer \cite{RaySinger}. Avoiding any
technical point at this stage, we simply recall that such a
definition is build from an analytic continuation of $ s \mapsto
\mbox{tr} (\Delta_g-z)^{s} $, using that $ (\Delta_g-z)^{s} $ is
trace class at least for $ \mbox{Re}(s) $ sufficiently negative,
which uses crucially  the discreteness of the spectrum of $
\Delta_g $.

In this spirit, the first goal of this paper is to realize the
{\it resonances} of Schr\"odinger
operators  with AC spectrum, as the zeroes of a determinant defined via a certain
Zeta function. 

Let us informally recall that, if  $ H = H_0 +
V $ with $ H_0 = - \Delta $ on $ \Ra^d $ and  $ V $ a perturbation tending to $ 0 $ at infinity, the
resonances are the natural discrete spectral datum of the problem.
They can be defined as the poles of some meromorphic continuation
of the resolvent of $ H $ and thus can  be considered as the analogues of
the eigenvalues for confining systems. Notice however that, apart
from possible real eigenvalues, resonances usually refer to
complex poles.

The problem of defining resonances as zeroes of determinants is
very natural and has already been considered by several authors,
in connection with the important question of their distribution
\cite{Zw89,Vodev,Froese1,Froese2,PeZw99,PeZw01,Simon,GuZw,BoJuPe,BoBrRa,Gu06}.
In these references, various determinants are used such as absolute
determinants or relative determinants, determinants of the
scattering matrices. In this paper we will basically study
relative determinants. The corresponding construction is fairly
well known in the {\it relatively trace class} situation, ie when
$ (H-z)^{-k} - (H_0-z)^{-k} $ is of trace class, that is when $V$
decays sufficiently fast at infinity and we refer to \cite{Muller}
for a nice review on this case. The main point in this paper is to
consider determinants for slowly decreasing perturbations of long
range type. We first recall some well known facts.

When $ V = V(x) $ is a potential (or possibly a first order
differential operator), a natural candidate for our purpose can
be the so called {\it perturbation determinant} (see
\cite{Yafaev}) defined by
\begin{eqnarray}
D_p (z) = D_{p} (H_0,H;z) := \mbox{Det}_p  \left(
(H-z)(H_0-z)^{-1} \right) = \mbox{Det}_p \left( I + V (H_0-z)^{-1}
\right), \label{defFredholm}
\end{eqnarray}
where $ \mbox{Det}_p $ is the Fredholm determinant of order $p$
which is defined as follows (see \cite{GoKr1,Yafaev} for more
details). Given a separable Hilbert space (here $ L^2 (\Ra^d) $),
one defines the Schatten class of order $ p \geq 1 $ as the space
$ \mathbf{S}_p $ of compact operators $ K $ whose singular
numbers\footnote{ie the spectrum of $|K|:= (K^* K)^{1/2} $} form a
sequence in $ l^p (\Na) $ (for $ p = \infty $, $
\mathbf{S}_{\infty} $ is the class of compact operators). The most
classical examples are $\mathbf{S}_1 $, the trace class, and
$\mathbf{S}_2 $, the Hilbert-Schmidt class. Then, if
 $ K \in  \mathbf{S}_p $, the spectrum of $ K $ is also in $ l^p
(\Na) $ and, if $ p $ is an integer, one sets
\begin{eqnarray}
 \mbox{Det}_p (I+K) :=  \prod_{k \geq 0} (1 + \lambda_k)
\exp \left( \sum_{j=1}^{p-1} \frac{(-1)^j}{j} \lambda_k^j
 \right)  , \qquad (\lambda_k)_{k \geq 0} = \mbox{spec}(K) ,
 \label{detFredholm}
\end{eqnarray}
where the product is convergent since the Weierstrass function on
the right hand side is $ 1 + {\mathcal O}(\lambda^p_k) $. If $ V $
tends to zero with  rate $ \rho > 0 $, ie
\begin{eqnarray}
 |V(x)| \leq C \scal{x}^{-\rho}, \label{conditiondecroissancevague}
\end{eqnarray}
  it is classical that
\begin{eqnarray}
 V (H_0 - z)^{-1} \in \mathbf{S}_p \ \ \ \mbox{if} \ \ \min
(2,\rho) >  d/p . \label{Schattencondition}
\end{eqnarray}
 For instance, in dimension $d=1$ with $ V $ of
short range, ie $ \rho >1 $, $ V (H_0 - z)^{-1} $ is  trace class
and one can define $ D_1 (H_0,H;z) $, which is essentially the
framework of \cite{Froese1,Simon}. The Fredholm determinant  of order $ 1
$ is a rather popular tool for several reasons. For instance, it
satisfies the formula
$$ \mbox{Det}_1 \left( (I + K_1)(I+K_2) \right) = \mbox{Det}_1 \left( I + K_1 \right)
\mbox{Det}_1 \left( I+K_2 \right) ,$$ as in finite dimension. This
formula doesn't hold for $ p \geq 2 $ (one needs then to add
correction factors). Also, formula (\ref{detFredholm}) shows that
for $ p = 1$, we have a 'pure' factorization of the determinant of
$ I + K $ by its eigenvalues $ 1 + \lambda_k $. It is nevertheless
necessary to consider Fredholm determinants of higher order.
Indeed, even for compactly supported potentials,  $ V (H_0 - z)^{-1} $
is not
 of trace class in general when $ d \geq 2 $ (basically $ V (H_0 - z)^{-k}  \in \mathbf{S}_1  $ if $ k > d / 2$
 and $ \rho > d $). Furthermore,
  even for $ d = 1 $, one also needs to consider $
 p \ne 1 $ to deal with long range potentials, ie when $ 0 < \rho \leq 1
 $.

\medskip

There is in addition a major drawback in the definition  (\ref{defFredholm}): it is restricted to relatively
compact perturbations. In particular, we can not consider $V$
that are second order differential operators.

One can overcome this difficulty by defining relative determinants
via relative Zeta functions. This construction was first
introduced for relatively trace class perturbations, ie basically
for perturbations with coefficients decaying like
(\ref{conditiondecroissancevague}) with $ \rho > d $ (see
\cite{Muller} for references) and was then
extended in \cite{BoucletAA,BoucletJFA} to general $ \rho
> 0 $, using an original idea of Koplienko \cite{Koplienko}. We
recall this construction. Let $ V $ be a differential operator of
the form
$$ V = \sum_{|\alpha| \leq 2} v_{\alpha} (x) D^{\alpha} , \qquad D = - i \partial_x , $$
symmetric on $ L^2 (\Ra^d) $ such that $ - \Delta + V $ is
uniformly elliptic, whose coefficients are smooth and satisfy
\begin{eqnarray}
 |\partial^{\beta} v_{\alpha}(x) | \leq C_{\beta} \scal{x}^{-
\rho}, \qquad x \in \Ra^d , \label{decroissanceprecise}
\end{eqnarray}
 for some $ \rho > 0 $. We shall
further on consider semiclassical operators, ie replace $ D $ by $
h D $ with $ h \in (0,1] $, and all the results quoted here for $
h = 1 $ will still hold. One defines  the so called  {\it regularized spectral shift
function} $ \xi_p \in {\mathcal S}^{\prime}(\Ra) $ (see \cite{BoucletAA,BoucletJFA})  as the unique
distribution vanishing near $ - \infty $ such that
\begin{eqnarray}
 \scal{\xi_p^{\prime},f} = \mbox{tr} \left( f (H_0 + V) -
\sum_{j=0}^{p-1} \frac{1}{j!} \frac{d^j}{d \varepsilon^j} f (H_0 +
\varepsilon V)_{| \varepsilon = 0 } \right),
\label{definitionscatteringphase}
\end{eqnarray}
for all  Schwartz function $ f  $, or more generally $ f \in
S^{-k}(\Ra) $ (ie $ \partial_{\lambda}^j f (\lambda) = {\mathcal
O}(\scal{\lambda}^{-k-j}) $) with $ k $ large enough. For $p=1$,
we recover the well known Kre$\check{\i}$n spectral shift
function. For $ p \geq 2 $, this trace regularization by Taylor's
formula is due to Koplienko \cite{Koplienko}. We also refer to the recent paper \cite{GPS} for a general introduction to Koplienko's regularized spectral shift function in connection with determinants. See also \cite{Koplienko2,Rybkin,Belov} in the one dimensional case.

Denoting by $ (\cdot - z)^{-s} $ the map $ \lambda \mapsto
(\lambda - z)^{-s} $, it is shown in \cite{BoucletJFA} that the
{\it regularized Zeta function},
$$ \zeta_p (s,z) := \scal{\xi^{\prime}_p, (\cdot-z)^{-s}}, \qquad \mbox{Im}(z)>0, \ \mbox{Re}(s)\gg 1 $$
has a meromorphic continuation, with respect to $s$, which is
regular at $s=0$. This allows to define
\begin{eqnarray}
D_p^{\zeta}(z) = D_p^{\zeta} \left(H_0, H_0 +V ; z \right) := \exp
\left( -
\partial_s \zeta_p (s,z)_{|s=0} \right) , \nonumber
% \label{detezeta}
\end{eqnarray}
 which is holomorphic for $ \mbox{Im}(z) > 0
$. The notation $ D_p^{\zeta} $ is justified by the fact that
\begin{eqnarray}
 D_p^{\zeta} \left(H_0, H_0 +V ; z \right) = D_p \left( H_0, H_0
+ V ; z \right), \label{coincidencedefinition}
\end{eqnarray}
 when $ V $ is a potential (see \cite{BoucletJFA}). In other words, the
definitions of the perturbation determinant by Fredholm
determinants and regularized Zeta functions coincide if they
both make sense. In addition, one proved in \cite{BoucletJFA}
that, in the distributions sense,
\begin{eqnarray}
 \frac{d}{d \lambda} \mbox{arg} \ D_p^{\zeta}(\lambda +
i \epsilon) \rightarrow - \pi \xi^{\prime}_p (\lambda), \qquad
\epsilon \downarrow 0 . \label{scatteringphase}
\end{eqnarray}
For this reason,  $ \xi_p $ is also called {\it generalized
scattering phase} of order $p$. The above formula is well known
for $ \xi_1 $ and was initially proved in \cite{Krei1} (see also \cite{BiKr}). See also
\cite{Koplienko,Koplienko2,GPS} for $ p \geq 2 $. Note the
parallel with the finite dimensional analogy of the very beginning
of this paper: for an hermitian matrix $ A $ on $ \Ca^N $ with
spectrum $ \lambda_1 , \ldots , \lambda_N $, one easily sees that
$$  \frac{d}{d \lambda} \  \mbox{arg} \ \mbox{Det}(A- \lambda - i \epsilon) \rightarrow  -
 \pi \sum_{ k = 1}^N \delta (\lambda - \lambda_k), \qquad \epsilon
\downarrow 0 , $$ where the right hand side is $ - \pi $ times the
derivative of the eigenvalue counting function, ie the analogue of
the spectral shift function for a discrete spectrum.
 This also suggests that if the resonances of $ H_0 +
V $ are indeed the zeroes of (a suitable meromorphic continuation
of) $ D_p^{\zeta}(z) $, the derivative of $ \xi_p
(\lambda) $ should involve a function (and/or a measure) with singularities carried by the resonances.
% (recall that the resonances may be complex but $ \xi_p $ is defined on $ \Ra $). 
Such a result is sometimes
referred to as Breit-Weigner formula and is already known for $ p
= 1 $ (see \cite{BruneauPetkov} and the references therein). In
this paper, we shall prove it for general $ p \geq 1 $. We will
also give semiclassical bounds.

\bigskip

Throughout this paper, we shall use the definition of resonances
and some related results given in \cite{Sjosequi} (see also
\cite{SjBook}). The definition is basically taken from  the
original paper by Sj\"ostrand-Zworski \cite{SjostrandZworski} and
the other useful results of \cite{Sjosequi} come from a
simplification of the proof of the trace formula \cite{Sjostrace}.
Before stating the results, we fix the notation and some definitions.

For $ 0 <  \theta_0 < \pi  $, $ R_0 > 0  $ and $ \epsilon_0
> 0   $, we set
\begin{eqnarray}
 \Sigma (\theta_0,R_0, \epsilon_0 ) := \{ r \omega \ ; \ \omega
\in \Ca^d, \ \mbox{dist}(\omega, {\mathbb S}^{d-1})<  \epsilon_0 ,
\ r \in e^{i[0,\theta_0]}(R_0,+\infty) \} . \nonumber
%\label{secteur}
\end{eqnarray}

\begin{defi} Let $ \rho > 0 $.
We define $ {\mathcal C}_{\rho}(\theta_0,R_0, \epsilon_0 ) $ as
the set of smooth functions   $  v $ on $ \Ra^d $ which have
 an analytic extension to $ \Sigma
(\theta_0,R_0, \epsilon_0 ) $ such that
\begin{eqnarray}
  |v (x)| \leq C \scal{x}^{-\rho}, \qquad x \in \Sigma
(\theta_0,R_0,\epsilon_0)  . \label{decroissanceanalytique}
\end{eqnarray}
Here $ \scal{x} = (1+|x|^2)^{1/2} $. A family  $
(v_{\iota})_{\iota \in I} $ is said to be bounded in $ {\mathcal
C}_{\rho}(\theta_0,R_0, \epsilon_0 ) $ if it is bounded in $
C^{\infty}(\Ra^d) $ and if the constant $ C $ in
(\ref{decroissanceanalytique}) is uniform with respect to $ \iota
\in I $.
\end{defi}

We  consider perturbations of $  H_0 (h) = - h^2 \Delta
 $ by second order differential operators of the form
\begin{eqnarray}
  V(h) = \sum_{|\alpha| \leq 2} v_{\alpha}(x,h) (h D )^{\alpha} , \label{coefficientoperateur} 
\end{eqnarray}
depending on a small parameter $ h > 0 $. We assume that, for some
$ h_0 > 0 $, the coefficients are such that, for all $ |\alpha|
\leq 2 $,
\begin{eqnarray}
   (v_{\alpha}(.,h))_{h \in (0,h_0]} \ \ \mbox{is bounded in} \ \
{\mathcal C}_{\rho}(\theta_0,R_0,\epsilon_0)  ,
\label{coeffapriori}
\end{eqnarray}
and such that, for some $ c > 0 $,
\begin{eqnarray}
 v_{\alpha}(.,h)  \mbox{ doesn't depend on}  \ h \  \mbox{if} \ |\alpha| = 2 ,
\label{symboleprincipal}  \\
  |\xi|^2 + \sum_{|\alpha| = 2} v_{\alpha}(x)\xi^{\alpha} \geq c
|\xi|^2, \qquad x \in \Ra^d, \ \xi \in \Ra^d . \label{ellipticite}
\end{eqnarray}
We also assume that
\begin{eqnarray}
 V(h) \ \mbox{is symmetric on } \ C_0^{\infty} (\Ra^d) \  . \label{autoadjoint}
\end{eqnarray}
These assumptions imply that $ H_0(h) + V(h) $ is selfadjoint on $
L^2 (\Ra^d) $ with domain $ H^2 (\Ra^d) $  the usual Sobolev
space.

The assumption (\ref{coeffapriori}) implies that the coefficients of $ V $ must be smooth on $ \Ra^d $. This is mostly for convenience, to simplify the analysis, but we expect that some local singularities could be considered as well, using for instance the black-box formalism of Sj\"ostrand-Zworski \cite{SjostrandZworski}. Notice however that, apart from the special case  $ p = 1 $, we have to consider operators of the form $ H_0 + \varepsilon V $ hence with $ H_0 $ and $ V $ defined on the same space. In particular, the generalization of the present results to perturbations by obstacles (+ long range metrics or potentials) would require a modified approach.

\medskip

\noindent {\bf Notation.} We shall mostly write $ H_0 $, $ V $ for
$ H_0 (h) $ and $ V (h) $.  When no confusion will be possible, $
V $ will also denote the family of operators $ (V(h))_{0 < h \leq
h_1} $. Such a family will sometimes be denoted by $ (V(h))_{ h
\ll 1} $ to mean that it is of the form $ (V(h))_{0 < h \leq h_1}
$ for some $ h_1 > 0$.

\medskip
 It
is convenient to summarize the above properties in the following
definition.
\begin{defi} We say that  $V = (V (h))_{h \in (0,h_1]} $
belongs to $ {\mathcal V}_{\rho}(\theta_0,R_0,\epsilon_0) $ if it
satisfies (\ref{coeffapriori}), (\ref{symboleprincipal}),
(\ref{ellipticite}) and (\ref{autoadjoint}). A family $
(V_{\iota})_{\iota \in I} = (V_{\iota}(h))_{h \in (0,h_1],  \iota
\in I } $ is bounded in $ {\mathcal V}_{\rho}
(\theta_0,R_0,\epsilon_0) $ if the families of coefficients $ (
v_{\alpha,\iota}(.,h) )_{h \in (0,h_1],  \iota \in I } $ are
bounded in $ {\mathcal C }_{\rho} (\theta_0,R_0,\epsilon_0) $ for all $ \alpha $ and if the
constant $ c $ in (\ref{ellipticite}) can be chosen independently
of $ \iota $.
\end{defi}

\noindent {\bf Remark.} To state this definition, we have
explicitly fixed the range of $ h $, namely $ (0,h_1 ] $, but we
will also freely write that $V = (V(h))_{h \ll 1} $ belongs to $
{\mathcal V}_{\rho}(\theta_0,R_0,\epsilon_0) $ to mean that, for
some $ h_1 $ small enough, $ (V(h))_{h \in (0,h_1]} \in {\mathcal
V}_{\rho}(\theta_0,R_0,\epsilon_0) $. A similar slight abuse of
notation will be used for families $  ( V_{\iota} )_{\iota \in I}
= ( V_{\iota}(h) )_{h \ll 1, \iota \in I} $.

\bigskip

Obviously,  any $ v \in {\mathcal
C}_{\rho}(\theta_0,R_0,\epsilon_0) $ satisfies
 (\ref{decroissanceprecise}). Therefore, using the
results of \cite{BoucletAA}, we can define the generalized
scattering phase $ \xi_p (.,h) $ associated to  $ -h^2 \Delta $
and $ V(h) $, provided
$$ p \rho > d . $$
We can then define the regularized Zeta function $ \zeta_p (s,z,h) $
by
$$  \zeta_p (s,z,h) = \scal{ \xi_p^{\prime}(.,h), (.-z)^{-s} }, \qquad \mbox{Im}(z)>0, \ \ \mbox{Re}(s) \gg 1 . $$
According to \cite{BoucletJFA}, $ \zeta_p (s,z,h) $ can be continued
analytically at $ s = 0 $ and we can define the relative
determinant of order $ p $
\begin{eqnarray}
 D_p^{\zeta} (z,h) := \exp \left( - \partial_s \zeta_p
(s,z,h)_{|s=0} \right), \qquad \mbox{Im}(z) > 0.  \label{detezeta}
\end{eqnarray}
We note that, for more precise purposes, the analytic continuation
(in $s$) of the Zeta function will be reviewed in Section
\ref{reviewzeta}.

The determinant $ D_p^{\zeta} (z,h) $  is our candidate to become
the 'characteristic polynomial' of the resonances of $ H_0 + V $.

 We now briefly recall
the definition of resonances of \cite{Sjosequi,SjostrandZworski}
(see Section \ref{OpRe} of the present paper for precise
statements).
  Let $ \theta_0 \in (0, \pi) $, $ \epsilon > 0 $  such that $ \epsilon < 2 \pi - 2 \theta_0
  $ and consider a relatively compact open subset
\begin{eqnarray}
 \Omega \Subset  e^{i (-2 \theta_0, \epsilon)}(0,+\infty)   \label{localiseOmega}
\end{eqnarray}
 which is
simply connected and such that
\begin{eqnarray}
\Omega \cap (0,+\infty)  \mbox{ is a non empty interval}.
\label{intervalle}
\end{eqnarray}
The resonances of $ H_0 + V $ in $ \Omega $ are by definition the
eigenvalues in $ e^{-i[0,2 \theta_0) (0,+\infty) } \cap \Omega $
of some non selfadjoint operator $ H_0 (\theta_0) + V (\theta_0) $
obtained by analytic distortion. We denote by
$$ \mbox{Res}(H_0+V,\Omega) := \mbox{set of resonances of} \
H_0 + V \  \mbox{ in } \  \Omega , $$ which is a finite set
depending on $h$. We recall here that, for the operators
considered in this paper, we have the following Weyl upper bound
for the number of resonances in $ \Omega $ (see for instance
\cite{Sjosequi}),
\begin{eqnarray}
\# \mbox{Res}(H_0+V,\Omega) \leq C h^{-d}, \qquad h \ll 1 .
\label{comptagedesresonances}
\end{eqnarray}
Note that they are counted with {\it multiplicity}  and that the
multiplicity of each resonance is  well defined as the rank of a
certain projector (see Section \ref{OpRe}) which is non orthogonal
in general.

%In the sequel, if $ U \subset \Ca $ is an open subset, we shall
%set
%$$ U^+ = U \cap \{ \mbox{Im}(z) > 0 \} , $$
%ie $ U^+ = U \cap \Ca^+ $.

\medskip

Our first result is the following.

\begin{theo} \label{theoremeprincipalp} Let $ \rho > 0 $, $ V \in {\mathcal
V}_{\rho} (\theta_0,R_0,\epsilon_0) $ and $ p > d / \rho $. Then,
for all  $h \ll 1 $, $ D_p^{\zeta}(z,h) $ has an analytic
continuation from
\begin{eqnarray}
 \Omega^+ := \Omega \cap e^{i(0,\epsilon)}(0,+\infty) \label{Omegasuperieur}
\end{eqnarray}
 to $ \Omega $, of the form
$$ D_p^{\zeta}(z,h) = \prod_{w \in {\rm Res}(H_0+V, \Omega)} (z - w) \times \exp (\varphi_p (z,h)) ,
\qquad z \in \Omega ,  $$  where each resonance is repeated
according to its multiplicity and the function $ z \mapsto
\varphi_p (z,h) $ is holomorphic on $ \Omega $.
\end{theo}

The proof is given in subsection \ref{preuvedutheoreme}.

Notice that the function $ \varphi_p (z,h) $ is uniquely defined
up to  a multiple of $ 2 i \pi $ of the form $ 2 i k(h)  \pi $. By
(\ref{scatteringphase}), an immediate consequence of Theorem
\ref{theoremeprincipalp} is the following Breit-Wigner formula.

\begin{coro} With the notation and assumptions of
Theorem \ref{theoremeprincipalp},  for all $ h \ll 1 $  we have
\begin{eqnarray}
\xi^{\prime}_p (\lambda,h) = \sum_{w \in {\rm Res}(H_0+V,\Omega)
\cap \Ra} \delta (\lambda-w) - \sum_{w \in {\rm Res}(H_0+V,
\Omega) \setminus \Ra} \frac{ {\rm Im}(w)}{\pi
|\lambda-w|^2} - \frac{1}{\pi} \emph{Im}(\partial_z
\varphi_p(\lambda,h)), \nonumber
\end{eqnarray}
  in the distributions sense on $ \Omega \cap
(0,+\infty) $.
\end{coro}

Here  $ \lambda $ is restricted to $ (0,+\infty) $, but it is well
known that
$$ \xi^{\prime}_p (\lambda,h ) = \sum_{w \in \sigma^- (H_0+V)} \delta (\lambda - w) , \qquad
  \lambda \in \Omega \cap (- \infty,0) , $$
where $ \sigma^-  (H_0 + V) = \sigma (H_0 + V) \cap (- \infty,0) $
is the set of negative eigenvalues of $ H_0 + V $ (see
\cite{BoucletAA} for instance but this is anyway an elementary
consequence of the definition (\ref{definitionscatteringphase})).

This corollary becomes of real interest if one has estimates on $
 \partial_z \varphi_p $. This is the purpose of the next
 results.

\begin{theo} \label{theoremeprincipal2} Assume that $ V \in {\mathcal
V}_{\rho} (\theta_0,R_0,\epsilon_0) $ with $ \rho > d / p $ and
$$ p = 1 \qquad \mbox{or} \qquad p = 2 . $$
Then any $ \varphi_p $ as in Theorem \ref{theoremeprincipalp}
satisfies, for any compact subset $ W \Subset \Omega $,
\begin{eqnarray}
| \partial_z \varphi_p (z,h) | \leq C_W h^{-d}, \qquad 
h \ll 1  , \ z \in W . \label{borneacontredire}
\end{eqnarray}
\end{theo}

This theorem is proved in subsection \ref{preuvedutheoreme2}. In Section \ref{sectionanalytique}, we also prove that a similar result holds for $ p \geq 3 $ if we assume that $ V $ is dilation analytic. However Theorem \ref{theoremeprincipal2} is
sharp in general for non globally analytic perturbations as is shown by  Theorem \ref{theoremeprincipal3} below.
% shows that
%(\ref{borneacontredire}) can not hold, in general, for $ p \geq 3$.

Fix first
$$ W = \{ z = r e^{- i \theta} \in \Ca \ ; \ 1 \leq r \leq 4 , \ 0 \leq \theta \leq \pi \} , $$
and observe that, for $ \pi/2 < \theta_0 < \pi $  and all $
\epsilon > 0 $ small enough, $ W $ is clearly contained in a
simply connected open set $ \Omega $ satisfying
(\ref{localiseOmega}) and (\ref{intervalle}). This neighborhood $
 \Omega$ can be chosen close enough to $ W $ so that we can define a determination of the square root $ z^{1/2} $,
  with  $ (r e^{-i\theta})^{1/2} = r^{1/2} e^{- i \theta/2} $ on $
W $ hence so that
$$ \mbox{Im}(z^{1/2}) \leq 0  \ \ \ \mbox{ on } \ \  W . $$
\begin{theo} \label{theoremeprincipal3} In dimension $d = 1 $, for all $ V \in  C_0^{\infty}(\Ra,\Ra) $, $ V \ne 0 $,
we can find $ \delta > 0 $ such that,
\begin{eqnarray}
 \limsup_{h \rightarrow 0} \sup_{z \in W} | h  e^{ \delta {\rm Im}(z^{1/2})/h} \partial_z
\varphi_3 (z,h) | = + \infty . \label{bornepresqueexponentielle}
\end{eqnarray}
In particular, $ | h   \partial_z
\varphi_3 (z,h) | $ can not be bounded on $ W $ uniformly with respect to $h$.
\end{theo}

The proof of this theorem is given in Section
\ref{preuvecontreexemple}.

 \medskip

We next give a general bound on $ \partial_z \varphi_p $ involving
the distorted operator $ H_0 (\theta) $ defined in Section
\ref{OpRe} and the semiclassical Sobolev space defined by
(\ref{Sobolevsemiclassique}). We recall that $ H_0 (\theta) -z$ is
invertible for all $ h \ll 1 $ and $ z \in \Omega $.

\begin{theo} \label{theoremeavecborne} Under the assumption of Theorem
\ref{theoremeprincipalp}, there exists $ N > 0 $ such that, for
all $ W   \Subset \Omega $,
$$ \left| \partial_z \varphi_p (z,h) \right|\leq C_W h^{-d} \sup_{Z \in \Omega}
\left( 1 + || (H_0(\theta_0)-Z)^{-1}||_{L^2 \rightarrow H^{2,0}_{\rm
sc }} \right)^{N} , \qquad h \ll 1, \ z \in W . $$
\end{theo}

 In general, $ || (H_0(\theta_0)-Z)^{-1}||_{L^2 \rightarrow
H^{2,0}_{\rm sc }} $ is
%expected to be
of order $ {\mathcal O}(e^{C h^{-1}}) $, locally uniformly with
respect to $Z$ (see Proposition \ref{propestires}). However,
Theorem \ref{theoremeprincipal2} shows that the corresponding
exponential upper bound on $ \partial_z \varphi_p $ can be much
improved if $ p = 1,2 $ (and $ p \geq 3 $ if $V $ is dilation analytic, see Section \ref{sectionanalytique}). Note also that Theorem
\ref{theoremeprincipal3} can be interpreted as a weak exponential
lower bound.

Theorem \ref{theoremeavecborne} is proved in subsection
\ref{preuvedutheoreme}.

%the following definition will be convenient.
%\begin{defi} The space $ {\mathcal H}_{\rm pol}(\Omega,H_0) $ (resp. $ {\mathcal H}_{\rm exp}(\Omega,H_0) $) is
%the space of $h$-dependent holomorphic functions on $ \Omega $
%such that, for all $ W \Subset \Omega $,
%$$ |g (z,h)| \leq C_{W} h^{-d}, \qquad z \in W, \ h \in (0,H_0 ], $$
%$$ \left( \mbox{resp.} \qquad  | g (z,h)| \leq e^{- C_{W} h^{-1} }, \qquad z \in W, \ h \in (0,H_0
%]. \right) $$
%\end{defi}

\bigskip

To motivate the analysis developed in the next sections, let us
already show that most of the results above will essentially be
reduced to the study of $ \zeta_p (k,z,h)  $, for some $k$ large
enough.

% As we shall see below in Proposition
%\ref{reductionZetaintro}, Theorems \ref{theoremeprincipalp} and
%\ref{theoremeprincipal2}

The basic strategy is the following. Using (\ref{detezeta}),  we have
\begin{equation}\label{dkdetzeta}
\partial^{k}_z \log   D_p^{\zeta} (z,h) = - \partial^{k}_z
\partial_s \zeta_p (s,z,h)_{|s=0},
\qquad k \geq 1 , \ \ z \in \Omega^+ .
\end{equation}
Here and below $ \partial_z^k \log g  $ stands for  $ \partial_z^{k-1} ( g^{\prime}/g ) $,  for any non
vanishing holomorphic function $g$.
On the other hand, at least for $k > d/2$, we also have
\begin{eqnarray} \label{derivezeta}
\partial^{k}_z  \partial_s \zeta_p (s,z,h)_{|s=0}=(k-1)!\,  \zeta_p (k,z,h),
\end{eqnarray}
as will be proved in Section \ref{reviewzeta} (see (\ref{derivezeronu}) and the discussion thereafter)
and is formally a consequence of the identity,
\begin{eqnarray}
\partial_z^{k} \partial_s (\lambda-z)^{-s}_{\vert s=0} =  (k-1)! (\lambda-z)^{-k}. \label{pourleformel}
\end{eqnarray}
Fix then $z_0 \in \Omega^+  $. In Section
\ref{reviewzeta} (see Proposition \ref{astucevincent}) we shall also prove that, for all  $ \nu \geq 0  $,
\begin{eqnarray}
 |  \partial_z^{\nu+1}   \partial_s \zeta_p (s,z_0,h)_{|s=0} |  \leq C h^{-d} , 
 \qquad h \ll 1 .
% \in (0,h_1] .
\label{upperzeta}
\end{eqnarray}
%Note that the case $ \nu = 0 $ is allowed in  (\ref{upperzeta}). 
In addition, by
(\ref{comptagedesresonances}), we have, for all $ \nu \geq 0 $,
\begin{eqnarray}
\sum_{w \in {\rm Res}(H_0+V,\Omega)} |z_0-w|^{-\nu-1} \leq C h^{-d},
\qquad h \ll 1
% \in (0,h_1] 
, \label{constanteWeyl}
\end{eqnarray}
since $ |z_0-w| \gtrsim 1   $. These are the
essential tools of the reduction given by Proposition
\ref{reductionZetaintro} below. Before stating it and to
 consider the different possible estimates for $ \partial_z
\varphi_p $, we introduce the following. Let
$$ {\mathcal H}_{\rm
hol} (\Omega,h_1) := \{  (\phi (.,h))_{h \in (0,h_1]}  \} , $$ be
the space  of $h$-dependent  families of holomorphic functions on
$ \Omega $. Let $ {\mathcal H} (\Omega,h_1) $ be a subspace of $
{\mathcal H}_{\rm hol} (\Omega,h_1) $ such that
\begin{eqnarray}
 (h^{-d})_{h \in (0,h_1]} \in {\mathcal H} (\Omega,h_1) ,
\label{stabiliteconstantes}
\end{eqnarray}
and which is stable by taking the primitive, ie such that for all $ (\phi (.,h))_{h \in (0,h_1]} \in {\mathcal H}_{\rm
hol} (\Omega,h_1) $ and some $ z_0 \in \Omega $,
\begin{eqnarray}
 (\phi^{\prime}(.,h))_{h \in (0,h_1]} \in {\mathcal H} (\Omega,h_1) \qquad 
\Rightarrow \qquad  (\phi (.,h) - \phi (z_0,h))_{h \in (0,h_1]} \in {\mathcal H} (\Omega,h_1). \label{stabiliteparprimitive}
\end{eqnarray}
Note that, if $ z_0 $ is such that $  |\phi (z_0,h) |
\lesssim h^{-d} $, and by  using (\ref{stabiliteconstantes}), one can replace (\ref{stabiliteparprimitive})  by $ (\phi^{\prime}(.,h))_{h \in (0,h_1]} \in {\mathcal H} (\Omega,h_1)  
\Rightarrow  (\phi (.,h) )_{h \in (0,h_1]} \in {\mathcal H} (\Omega,h_1) $.

\medskip 

\noindent {\bf Example.} The space  $
{\mathcal H}_{\rm hol} (\Omega,h_1) $ itself or the subspace of
functions such that, for all $ W \Subset \Omega $, $ |\phi (z,h) |
\leq C_W h^{-d} $ for all $ z \in W $ and $ h \in (0,h_1 ] $ satisfy (\ref{stabiliteconstantes}) and (\ref{stabiliteparprimitive}).

 %The latter is a simple consequence of the Cauchy formula. 

\begin{prop} \label{reductionZetaintro} If we can find $ h_1 > 0 $ small enough,  $k \geq 1$ and $
\phi_p \in {\mathcal H} ( \Omega , h_1 ) $ such that
\begin{eqnarray}
\zeta_p (k,z,h) = \sum_{w \in {\rm Res}(H_0+V,\Omega)}
\frac{1}{(w-z)^k} + \phi_p(z,h) , \qquad z \in \Omega^+  , \ h \in
(0,h_1 ],
  \label{aintegrerkfois}
\end{eqnarray}
then Theorem \ref{theoremeprincipalp} holds true with $ \varphi_p
$ such that $ \partial_z \varphi_p \in  {\mathcal H} ( \Omega ,
h_1 ) $.
\end{prop}

\medskip

\noindent {\it Proof.} Setting for simplicity
$$ D  = D_p^{\zeta}
\left( z ,h \right) , \qquad F  = \prod_{w \in
    \rm{Res}(H_0+V,\Omega)} (z-w) , $$
which are holomorphic and don't vanish on  $ \Omega^+  $,
 (\ref{dkdetzeta}), (\ref{derivezeta}) and  (\ref{aintegrerkfois})
give
\begin{eqnarray}
  \partial_z^{k-1} \left( \frac{\partial_z D }{D} -
\frac{\partial_z F }{F} \right) = -(k-1)! \phi_p , \qquad \mbox{on} \ \Omega^+
.  \label{integrablefacile}
\end{eqnarray}
If $ k = 1 $, we therefore obtain 
\begin{eqnarray}
   \frac{\partial_z D }{D} -
\frac{\partial_z F }{F} \in {\mathcal H} (\Omega , h_1) , \label{sousconclusion}
\end{eqnarray}
 which implies easily  the result. If $ k-1 \geq 1$, we denote by $ \Phi_p $ the $( k-1)$-th primitive
of $ - (k-1)! \phi_p  $ (ie
  $\partial_z^{k-1} \Phi_p =  - (k-1)! \phi_p$) such that
$$ \partial_z^{\nu} \Phi_p (z_0,h) =  \partial_z^{\nu} \left( \frac{\partial_z D }{D} -
\frac{\partial_z F }{F} \right) (z_0,h)  , \qquad 0 \leq \nu \leq
k-2 ,  $$
where $ z_0 $ is chosen arbitrarily in $ \Omega^+ $. The existence and uniqueness of $ \Phi_p  $ is
guaranteed by the simple connectedness of $ \Omega  $. By (\ref{upperzeta}) and (\ref{constanteWeyl}), we have
\begin{eqnarray}
| \partial_z^{\nu} \Phi_p (z_0,h) | \leq C h^{-d} , \nonumber
\end{eqnarray}
and this implies, together with (\ref{stabiliteconstantes}) and (\ref{stabiliteparprimitive}),  that
$$ \phi_p \in {\mathcal H}(\Omega, h_1) \qquad \Rightarrow \qquad \Phi_p \in {\mathcal H}(\Omega , h_1) .  $$
Thus (\ref{integrablefacile}) imply that (\ref{sousconclusion}) holds also if $ k - 1 \geq 1 $ and we get the result. \finpreuve

\section{The Zeta function}  \label{reviewzeta}
\setcounter{equation}{0} In this subsection, we  review the
construction of the meromorphic continuation of $ s \mapsto
\zeta_p (s,z,h) $. Although the latter was shown in
\cite{BoucletJFA} (for fixed $h$), we need to review the main
lines of the proof in order to prove the identity
(\ref{derivezeta}) and the estimate (\ref{upperzeta}).

We start with general considerations. Using the principal
determination of $ \log $ on $ \Ca \setminus (-\infty,0] $, we can
define $ (\lambda - z )^{-s} $ for $ s \in \Ca $, $ \lambda \in
\Ra $ and
 $ z \in \Ca \setminus [ \lambda, + \infty ) $. One can
then check that
\begin{eqnarray}
(\lambda-z)^{-s} = \frac{1}{\Gamma(s)} \int_0^{+\infty} t^{s-1}
e^{-t(\lambda-z)} d t, \qquad \mbox{Re}(z)<\lambda, \
\mbox{Re}(s)>0 ,
\end{eqnarray}
since both sides are holomorphic with respect to $z$ and the
equality holds for $ z \in (-\infty,\lambda) $ by an elementary
change of variables in the definition of $\Gamma(s)$. Next, if $ u \in {\mathcal S}^{\prime}(\Ra) $
is a temperate distribution such that, for some $ \lambda_0 > 0 $,
\begin{eqnarray}
 \mbox{supp}(u) \subset [ \lambda_0 , + \infty )
\end{eqnarray}
we can consider its Laplace transform $ L u (t) :=
\scal{u,e^{-t.}} $ ($ e^{-t.} $ stands for the map $
 \lambda \mapsto e^{-t \lambda}$),
%\footnote{where $ e^{-t.} $ stands for the map $ \lambda \mapsto e^{-t \lambda}$},
 and, for all $ \delta
> 0 $,
\begin{eqnarray}
| L u (t) | \leq C_{\delta} e^{-t (\lambda_0 - \delta)} , \qquad t
> 0 .
\end{eqnarray}
Furthermore, using that $ | \scal{u,f} | \leq C \sup_{j+k \leq N}
\sup_{\lambda \in \Ra} |\lambda^j \partial_{\lambda}^k f
(\lambda)| $ for some $ N $ and all $ f \in {\mathcal S}(\Ra) $, $
\scal{u,f} $ is still well defined if $ f (\lambda) = (\lambda -
z)^{-s} $ with $ \mbox{Re}(s) > s_0 $ large enough and $
\mbox{Re}(z) < \lambda_0 $. If in addition, we know that
\begin{eqnarray}
| L u (t) | \leq C t^{-d/2} , \qquad t \in (0,1]
\label{controleenzero}
\end{eqnarray}
then, one has
$$ \scal{u,(\cdot - z)^{-s}} =  \frac{1}{\Gamma(s)} \int_0^{+\infty} L u (t) e^{t z} t^{s-1}
 d t, \qquad \mbox{Re}(z)<\lambda_0 , \ \mbox{Re}(s)> \max
( s_0 , d /2 ) . $$ Note that the power $d/2$ could actually be
any arbitrary real number but, in the applications below, we shall
need only to consider this case. If (\ref{controleenzero}) is
replaced by the stronger assumption that there is an asymptotic
expansion at $t=0$, namely that, for all $ J > 0 $,
\begin{eqnarray}
Lu(t)= \sum_{j < J} a_j t^{-d/2 + j/2} + t^{-d/2 + J/2} b_J (t),
\qquad |b_J (t)| \leq C , \ \ t \in (0,1] ,
\label{developpementchaleur}
\end{eqnarray}
then we can write, for $\mbox{Re}(z)<\lambda_0$ and  $ \mbox{Re}(s)>
\max ( s_0 , d /2 ) $,
\begin{eqnarray}
\scal{u,(\cdot - z)^{-s}} =   I (s,z) +
II_J (s,z) + III_J (s,z) ,
\end{eqnarray}
with
\begin{eqnarray}
I (s,z) & = & \frac{1}{\Gamma(s)} \int_1^{\infty} L u (t) e^{tz} t^{s-1} dt , \nonumber \\
II_J (s,z) & = & \frac{1}{\Gamma(s)}  \int_0^1 b_J (t) e^{tz} t^{-d/2+J/2+s-1} dt , \nonumber \\
III_J (s,z) & = & \frac{1}{\Gamma(s)} \sum_{j < J} a_j \int_0^1 e^{tz} t^{-d/2 + j/2 +
s -1} dt . \nonumber
\end{eqnarray}
By choosing $ J > d $, both $ I $ and $ II_J $ are holomorphic
close to $ s = 0 $. Thus, using the fact that $ d \Gamma^{-1}(s) /
ds = 1 $ at $ s = 0$ and that $ \Gamma^{-1}(s) $ vanishes
at $0$  one sees that, for all $k \geq 1 $,
\begin{eqnarray}
\partial_z^{k} \partial_s  F (s,z)_{| s=0}= \Gamma (k) F (k , z) =  (k - 1)! F (k,z), \qquad \mbox{Re}(z)<
\lambda_0, \label{FIetII}
\end{eqnarray}
for $ F = I$ and $ F = II_J $. The term $ III_J $ can be computed
explicitly, namely,
\begin{eqnarray}
\Gamma (s) \times III_J (s,z) = \sum_{j = 0}^{ J-1} a_j \sum_{l \geq 0}
\frac{z^l}{l!} \frac{1}{s+j/2+l-d/2} . \label{serie}
\end{eqnarray}
At $s = 0$,
there is at most a simple pole, which corresponds to the terms
where $ j/2 + l - d/2 = 0$. Thus $ III_J (s,z) $ is
regular at $s=0$. This shows the existence of a meromorphic
continuation to the complex plane for
$$ s \mapsto \scal{u,(\cdot - z)^{-s}} = : Z (s,z) , $$
which is regular at $ s = 0 $.
 Furthermore one has,
\begin{eqnarray}
\partial_z^{k} \partial_s  III_J (s,z)_{| s=0}= (k - 1)! III_J (k,z), \qquad  k >  d/2,
\end{eqnarray}
(with $ k $ integer) since this derivative only involves terms with $ l > d/2 $ in
(\ref{serie}). Hence, using (\ref{FIetII}), we also have
\begin{eqnarray}
 \partial_z^{k} \partial_s Z (s,z)_{|s=0} = (k - 1)! Z
(k,z) , \qquad \mbox{Re}(z) < \lambda_0 , \ \ k >
d/2 . \label{derivezeronu}
\end{eqnarray}
Note that, if $ u $ is compactly supported, (\ref{derivezeronu})
is a direct  consequence of the identity (\ref{pourleformel}).

\medskip

When $ u = \xi_p^{\prime} $, the existence of a meromorphic
continuation in $s$ for $ \zeta_p (s,z,h) $ is a
consequence of the existence of an expansion of the form
 (\ref{developpementchaleur}) proved in \cite{BoucletAA}. Notice that altering $ L u (t) $ by an analytic function in $t$
will not destroy the form of this expansion. There is no
 restriction on $ \mbox{Re}(z) $ since, for all $ \lambda_0 \in \Ra $,
  $ \xi^{\prime}_p $ can be written as the sum of a compactly
  supported distribution and a temperate distribution supported in
  $ [ \lambda_0 , + \infty ) $ for which
  (\ref{developpementchaleur}) still holds since the Laplace
  transform of the compactly supported distribution is analytic in
  $t$.

In particular, for $ u = \xi_p^{\prime} $, the relation (\ref{derivezeronu}) yields (\ref{derivezeta}).
%\begin{equation} %\label{Zzeta}
%\partial^{k}_z  \partial_s \zeta_s (z,h)_{|s=0}=(k-1)!\,  \zeta_{k}(z,h), \qquad k >\frac{d}{2} \nonumber
%\end{equation}
%which is exactly ().

\medskip

We now consider (\ref{upperzeta}).

\begin{prop} \label{astucevincent} For all $ z_0 \in \Omega^+ $ and all integer $ \nu \geq 0
$, (\ref{upperzeta}) holds.
\end{prop}

\noindent {\it Proof.}  We shall see that the result follows
from the following two facts: the existence of a semi-norm  $ || .
||_{\mathcal S} $  (independent of $h$) of the Schwartz space $
{\mathcal S}(\Ra) $ such that
\begin{eqnarray}
 | \left\langle \xi^{\prime}_p(h), \psi \right\rangle | \leq C h^{-d}|| \psi ||_{\mathcal
 S}, \qquad \psi \in {\mathcal S}(\Ra), \ h \in (0,h_0] , \label{equicontinuite}
\end{eqnarray}
and the existence of an expansion of the form
\begin{eqnarray}
\left\langle \xi^{\prime}_p(h),e^{-t ( . )} \right\rangle \sim
t^{-d/2} \sum_{j \geq 0} a_j (h) t^{j/2}, \qquad t \rightarrow 0,
\qquad \mbox{with} \ \  a_j (h) = {\mathcal O}(h^{-d}) .
\label{doubleexpansion}
\end{eqnarray}
The latter means that the difference between the left hand side
and the sum truncated at the order $M$ is bounded by $ C h^{-d}
t^{(M-d)/2} $, for $ t \in (0,1] $ and $ h \in (0,h_0]$. Indeed,
by writing $ \xi^{\prime}_p = \chi \xi^{\prime}_p + (1 - \chi)
\xi^{\prime}_p $
 with $ \chi \in C_0^{\infty}(\Ra) $ such that $ \chi \equiv 1 $ on a large enough compact set,
we may assume that $ (1-\chi) \xi^{\prime}_p $ is supported in $ [
\lambda_0 , + \infty ) $ with $ \lambda_0 > \mbox{Re}(z_0) $.
Therefore, using (\ref{equicontinuite}), (\ref{doubleexpansion}) and the discussion prior to Proposition \ref{astucevincent}, we see that $ \scal{\chi \xi^{\prime}_p(h), (\cdot-z)^{-s}} $ as well
as the terms $ I(h),II_J(h),III_J(h) $ corresponding to $ u = u(h) = (1 - \chi)
\xi_p^{\prime}(h) $ are $ {\mathcal O}(h^{-d}) $ uniformly with
respect to $ s $ close to $ 0 $ and $ z $ close to $z_0$ which
gives the result.

The proof of (\ref{equicontinuite}) can be found in
\cite{BoucletAA} so we only consider
(\ref{doubleexpansion}).
  For the latter, the main
remark is that, for all $ \varepsilon \in [0,1] $,
$$ - t (H_0 + \varepsilon V) =  (h t^{1/2} )^2 \Delta - \varepsilon \tilde{V}(h,t^{1/2},x, h
t^{1/2}D) $$ with
$$ \tilde{V}(h,t^{1/2},x,\xi) =  \sum_{l=0}^2
t^{1-\frac{l}{2}} \sum_{|\alpha|=l} v_{\alpha}(x,h) \xi^{\alpha} $$ where
the $ v_{\alpha}$ are defined by (\ref{coefficientoperateur}). By reviewing the
proof of Proposition 3.1 in \cite{BoucletAA} with $ h t^{1/2} $ as
 new semi-classical parameter, we see that, for all $M$, we have
the following expansion
\begin{eqnarray*}
 \mbox{tr} \left( e^{-t (H_0 + V)} - \sum_{j=0}^{p-1}
\frac{1}{j!} \frac{d^j}{d \varepsilon^j} e^{-t (H_{0}+\varepsilon
V)} |_{\varepsilon =0} \right) & = & \sum_{q <
M} (h t^{1/2})^{q-d} d_q (t^{1/2},h)  \\
& &  \qquad \qquad + (h t^{1/2})^{M-d} R_M (t^{1/2},h) ,
\end{eqnarray*}
with $ R_M (t^{1/2},h) = {\mathcal O}(1) $ for $ h \in (0,h_0] $
and $ 0 < t \leq 1 $. The coefficients $d_q (t^{1/2},h)$ are
smooth at $ 0 $ with respect to $ t^{1/2} $ and bounded with
respect to $h \in (0,h_0]$ as well as their derivatives so $
(\ref{doubleexpansion} $) follows. \finpreuve

\section{Trace class estimates}
\setcounter{equation}{0} In the sequel, we shall use the notation
$ O \! p_h^w (a) $ for standard $h$-pseudodifferential operators
of the form
$$ O \! p_h^w (a) u (x) = (2 \pi)^{-d} \int \! \! \int e^{i (x-y) \cdot \xi} a \left( \frac{x+y}{2} ,h \xi \right)
u(y) d\xi d y , \qquad h \in (0,h_0] , $$ with symbols $ a \in
S^{\mu,\nu} $, $ \mu,\nu \in \Ra $, namely such that
$$ |
\partial_x^{\alpha}
\partial_{\xi}^{\beta} a (x,\xi) | \leq C_{\alpha \beta} \scal{x}^{\nu}
\scal{\xi}^{\mu-|\beta|}. $$ We refer for instance to
\cite{RoBook,Mart,DiSj} for the proofs of the standard results we
shall use below on the analysis of such operators.
 We equip $S^{\mu,\nu} $ with its
standard Fr\'echet space topology given by the seminorms defined
by the best constants $ C_{\alpha \beta} $.

 We also define the
following semiclassical weighted Sobolev spaces
$$ H_{\rm sc}^{s,\sigma} := \scal{x}^{- \sigma} \scal{h D}^{-s} L^2 (\Ra^d), \qquad  s, \sigma \in \Ra , $$
equipped with the $h$-dependent norm
\begin{eqnarray}
 ||u||_{H^{s,\sigma}_{\rm sc}} := || \scal{hD}^s
\scal{x}^{\sigma} u ||_{L^2(\Ra^d)} . \label{Sobolevsemiclassique}
\end{eqnarray}
 Notice that
$$ H_{\rm sc}^{s,\sigma} \subset H_{\rm sc}^{s,0} \subset L^2 (\Ra^d), \qquad \mbox{if} \ s \geq 0, \ \sigma \geq 0 . $$
In this section,  we will consider  $h$-dependent families of
symbols
$$ a = (a(h))_{h \in (0,h_0]}, \qquad a (h) \in S^{2,0} \ \mbox{for all} \ \ h \in (0, h_0 ] . $$
Most of the time, we shall assume the existence of $C > 0$ such that, for
all $ h \in (0,h_0 ] $,
\begin{eqnarray}  |a(h,x,\xi)| \geq
C^{-1} |\xi|^2  , \qquad  x \in \Ra^d , \  |\xi| > C .
\label{elliptiquePDE}
\end{eqnarray}
When $ a = (a(h))_{h \in (0,h_0]} $ or $ b = (b(h))_{h \in (0,h_0]}
$, we shall adopt the short notation
$$ A = O \! p_h^w (a(h)), \qquad B = O \! p_h^w (b(h)) ,   $$ for all $ h \in (0, h_0 ] $. 
%This uniform condition will only be required for the estimate
%(\ref{normrestriction}) below.

In the next proposition, $ {\mathcal B} $ denotes a subset of $ (
S^{2,0})^{(0,h_0]} $, namely a set of families $ (a(h))_{h \in (0,h_0]} $, uniformly
bounded in $ S^{2,0} $, ie such that $ \{ a (h) \ ; \ h \in (0,h_0 ] ,
a \in {\mathcal B} \} $ is bounded in $ S^{2,0} $. We also assume that (\ref{elliptiquePDE})
holds for all $ a \in {\mathcal B} $, with a constant $ C > 0 $ independent of $a$.

%For instance, if $ (a(h))_{h \in (0,h_0]} $ is a bounded family in
%$ S^{2,0} $ satisfying (\ref{elliptiquePDE}), then for all bounded
%subset $ U \subset \Ca $, the set $ {\mathcal B} = \{ (a(h) -
%z)_{h \in (0,h_0] } \ ; \  z \in U \} $ satisfies this
%condition too, with perhaps a larger $C$.

%If, for all $ h \in (0,H_0 ]$, we have a sequence $ (
%a_n (h) )_{n \geq 0} $ converging to $ a(h) $, then under
%reasonable, conditions on these sequences, the set $ {\mathcal
%A}_2 = \{ a_n = (a_n(h))_{h \in (0,h_0]} \ ; \ n \geq 0 \} \cup \{
%a \} $ will also work.

%For $ m \leq 2 $ and $ \mu \leq 0 $, we will also consider sets $
%{\mathcal B}^{m,\mu} $ of families $ b = (b (h))_{h \in (0,H_0 ]}
%$ such that
%$$ \{ b_h \ ; \ h \in (0,h_0] , \ b \in {\mathcal B}^{m , \mu}  \} \subset S^{m,\mu} \ \mbox{is bounded} . $$
%For any $ b \in {\mathcal B}^{m,\mu} $ and $ h \in (0,H_0 ] $, we
%set

\begin{prop} \label{inverse}
Assume that, for all $ a \in {\mathcal B} $ and all $ h \in (0,h_0
] $,
$$ A  : H^{2,0}_{\rm sc} \rightarrow L^2
(\Ra^d) \ \mbox{is invertible} . $$ Then, for all $ s \geq 0 $ and
$ \sigma \geq 0 $, the restriction $$ A_{s,\sigma} = A |_{
H^{s+2,\sigma}_{\rm sc} }
$$
 is bounded from $ H^{s+2,\sigma}_{\rm sc} $ to $ H^{s,\sigma}_{\rm sc} $
with bounded inverse such that
\begin{eqnarray}
  A_{s,\sigma}^{-1} = A^{-1} |_{H^{s,\sigma}_{\rm sc}} . \label{restrictioninverse}
\end{eqnarray}
Furthermore, there exists $ C_{s,\sigma} > 0 $ such that, for all
$ h \in (0,h_0] $ and all $ a \in {\mathcal B} $,
\begin{eqnarray}
|| A_{s,\sigma}^{-1} ||_{H^{s,\sigma}_{\rm sc} \rightarrow
H^{s+2,\sigma}_{\rm sc}} \leq C_{s,\sigma}
% || A^{-1} ||_{L^2 \rightarrow H^{2,0}_{\rm sc}} 
\left( 1 + || A^{-1} ||_{L^2
\rightarrow H^{2,0}_{\rm sc}} \right)^{[\sigma]+1} ,
\label{normrestriction}
\end{eqnarray}
with $ [\sigma] $ the smallest integer $ \geq \sigma $.
\end{prop}
The equality (\ref{restrictioninverse}) means that we can consider
$ A^{-1} $ as an operator from $ H^{s,\sigma}_{\rm sc} $ into $
H^{s+2,\sigma}_{\rm sc} $ and (\ref{normrestriction}) gives an
estimate on the corresponding norm. Abusing the notation, this
proposition will allow us to denote  $ A^{-1} $ instead of $
A^{-1}_{s,\sigma} $ in the sequel.

\bigskip

\noindent {\it Proof.} The boundedness of $ A_{s,\sigma} $ follows
 from the $ L^2 $ boundedness of
$$ \scal{h D}^{s} \scal{x}^{\sigma} O \! p_h^w (a (h)) \scal{x}^{-\sigma} \scal{h D}^{-s-2} =: O \! p_h^w (b_{s,\sigma}(h)) $$
since $ b_{s,\sigma}(h) $ so defined belongs to $ S^{0,0} $. If $
\sigma > 0 $, we  consider next $ \sigma_1 := \sigma / [\sigma]
\in [0,1] $. Then, by the resolvent identity,
$$ A^{-1} \scal{x}^{\sigma_1} = \scal{x}^{\sigma_1} A^{-1} - A^{-1} [A,\scal{x}^{\sigma_1}]  A^{-1}  $$
where $ [A,\scal{x}^{\sigma_1}] = O \!  p_h^w (a_{\sigma_1}(h)) $
for some symbol $ a_{\sigma_1}(h) \in S^{1,0} $ depending
continuously on $a(h)$. Thus
$$ \scal{x}^{-\sigma_1}  A^{-1}  \left( 1 + [A,\scal{x}^{\sigma_1}]  A^{-1} \scal{x}^{-\sigma_1} \right) = A^{-1}  \scal{x}^{-\sigma_1} $$
 shows that $ A^{-1} $ is bounded from $ H^{0,\sigma_1} $ to $ H^{0,\sigma_1}
 $ with norm controlled, uniformly  with respect to $ a \in {\mathcal B} $ and $ h \in (0,h_0 ] $, by
 $ ||A^{-1} ||_{L^2 \rightarrow H^{2,0}_{\rm sc}} (1+||A^{-1}||_{L^2 \rightarrow H^{2,0}_{\rm sc}})
 $. By iteration, we obtain that $ A^{-1} $ maps continuously $ H^{0,2\sigma_1}_{\rm
 sc}, H^{0,3\sigma_1}_{\rm
 sc}, \ldots , H^{0,[\sigma]\sigma_1}_{\rm
 sc}
 $ into themselves and that
\begin{eqnarray}
  ||A^{-1} ||_{ H^{0,\sigma}_{\rm sc} \rightarrow H^{0,\sigma}_{\rm sc} } \leq C
 ||A^{-1} ||_{L^2 \rightarrow H^{2,0}_{\rm sc}} (1+||A^{-1}||_{L^2 \rightarrow H^{2,0}_{\rm
 sc}})^{[\sigma]}, \label{apriorisigma}
\end{eqnarray}
with $ C $ independent of $h$ and of $ a \in {\mathcal B} $. Using (\ref{elliptiquePDE}),  we
can construct, for all $ N \geq 0 $, symbols $ \tilde{a}_N(h) \in
S^{-2,0} $ and $ r_N (h) \in S^{-N,0} $, depending continuously on
$ a (h) $, such that
$$ O \! p_h^w (\tilde{a}_N(h)) O \! p_h^w (a(h))
= 1 + O \! p_h^w (r_N(h)) .$$ Notice that this is not a
semiclassical parametrix (that would be the case if we had a
remainder of the form $ h^N O \! p_h^w (r_N(h)) $) since
(\ref{elliptiquePDE}) is not an ellipticity condition in the
semiclassical sense. This is simply an $h$-dependent
classical parametrix (in the sense of Theorem 18.1.9 of
\cite{Horm3}). The symbol  $ \tilde{a}_N (h) $ is constructed by
successive approximations starting from  $ (1-\chi) (\xi) / a
(x,\xi,h) $, with $ \chi \in C_0^{\infty} $ such that $ \chi (\xi)
= 1 $ for $ |\xi| \leq C $, and then following the usual iterative
scheme.  We then obtain
\begin{eqnarray}
 A^{-1} = O \! p_h^w (\tilde{a}_N(h))  -  O \! p_h^w (r_N(h)) A^{-1}
. \label{parametrixe1}
\end{eqnarray}
Since $ O \! p_h^w (\tilde{a}_N(h)) $ maps $ H^{s,\sigma}_{\rm sc}
$ into $ H^{s+2,\sigma}_{\rm sc} $   and $ O \! p_h^w (r_N(h)) $
maps $ H^{0,\sigma}_{\rm sc} $ into $ H^{N,\sigma}_{\rm sc} $ for all $ N \geq 0 $,
with norms uniformly bounded with respect to $ a $ and $ h $, the
right hand side of (\ref{parametrixe1}) is therefore bounded from
$ H^{s,\sigma}_{\rm sc } $ to $ H^{s+2,\sigma}_{\rm sc} $, by
choosing $ N \geq s + 2 $ and using (\ref{apriorisigma}). The
result then follows easily. \finpreuve

\bigskip

In the sequel we shall denote by $ {\mathcal L}({\mathcal
H}_1,{\mathcal H}_2) $ the Banach space of linear continuous
mapping between Hilbert spaces $ {\mathcal H}_1 $ and $ {\mathcal
H}_2 $, equipped with the usual norm. We also denote by $ {\mathcal L}_{\rm invertible}({\mathcal
H}_1,{\mathcal H}_2) $ the open subset of invertible mappings.

\begin{prop} \label{lemmeprop} Let $ a = (a (h))_{h \in (0,h_0]} $ be a family of $
S^{2,0} $ satisfying (\ref{elliptiquePDE})
 and let  $ U \subset \Ca  $ be an open subset.
Assume that
$$ A - z : H^{2,0}_{\rm sc} \rightarrow L^2 (\Ra^d) \ \ \
\mbox{is invertible} $$ for all $ z \in U $   and all $ h \in
(0,h_0] $.
%Then,

%\noindent i) for all $ s \geq 0
%$ and $ \sigma \geq 0 $,
%$$ U \ni z \mapsto (A-z)^{-1} \in {\mathcal L}(H^{s,\sigma}_{\rm sc},H^{s+2,\sigma}_{\rm sc})
%$$
%is holomorphic.

\noindent i) Let $ b = (b(h))_{h \in (0,h_0]} $ be a family of $
S^{2,0} $. Then, for all $ h \in (0,h_0] $ and all $ z_0 \in U $, there exists $ \varepsilon_{h,z_0}
> 0 $ and a neighborhood $ U (z_0) \subset U $ of $z_0 $  such that, for all $ s , \sigma \geq 0 $, the map
\begin{eqnarray}
 (-  \varepsilon_{h,z} ,  \varepsilon_{h,z} ) \times U (z_0 ) \ni  (\varepsilon,z) \mapsto (A
+ \varepsilon B - z)^{-1} \in {\mathcal L}(H^{s,\sigma}_{\rm
sc},H^{s+2,\sigma}_{\rm sc}) \label{maplocal}
\end{eqnarray}
is well defined and smooth. In addition
\begin{eqnarray}
\frac{d^j}{d \varepsilon^j} (A + \varepsilon B-z)^{-1} = (-1)^j j
! (A + \varepsilon B-z)^{-1} \left( B (A + \varepsilon B-z)^{-1}
\right)^j . \label{deriveeTaylor}
\end{eqnarray}
ii) Assume that, for all $ h \in (0,h_0] $, we have a sequence $
( a_n (h) )_{n \in \Na}  $ converging to $  a (h) $ in $ S^{2,0}
$. Then, for all $ h \in (0,h_0] $ and all relatively compact
subset $ U_0 \Subset U $,  there exists $ n_{h,U_0} \in \Na $ such
that,
\begin{eqnarray}
A_n - z : H^{2,0}_{\rm sc} \rightarrow L^2(\Ra^d), \qquad z \in
U_0 , \ n \geq n_{h,U_0},
\end{eqnarray}
is invertible, and, for all $ s,\sigma \geq 0  $,
\begin{eqnarray}
 || (A_n - z)^{-1} - (A-z)^{-1} ||_{H^{s,\sigma}_{\rm sc}
\rightarrow
  H^{s+2,\sigma}_{\rm sc}} \rightarrow 0 , \qquad n \rightarrow
  \infty , \label{convergenceabstraite}
\end{eqnarray}
uniformly on $ U_0  $.
\end{prop}

\bigskip

\noindent {\it Proof.}  Fix $ h \in (0, h_0] $. Since $ B $ is
bounded from $ H^{2,0}_{\rm sc} $ to $ L^2 (\Ra^d) $,  for $
\varepsilon $ small enough and $ z $ close enough to $z_0$, $ A +
\varepsilon B - z $ is invertible. It is  then also invertible as a
bounded operator from $ H^{s+2,\sigma}_{\rm sc} $ to $
H^{s,\sigma}_{\rm sc} $ by Proposition \ref{inverse}. Since the
map $ T \mapsto T^{-1} $ is $ C^1 $ from $ {\mathcal L}_{\rm
invertible}(H^{s+2,\sigma}_{\rm sc}, H^{s,\sigma}_{\rm sc}) $ to $
{\mathcal L}(H^{s,\sigma}_{\rm sc}, H^{s+2,\sigma}_{\rm sc}) $,
(\ref{maplocal}) is $ C^1 $ with derivative given by
(\ref{deriveeTaylor}) with $ j = 1 $. The result then follows by
induction. Let us now prove {\it ii)}. Let $ z_0 \in U $. Since
 $ A - z_0 $ is invertible and by convergence of $ A_n $ to $ A
$, there exists $ n_{h,z_0} > 0 $ and $ \delta_{z_0,h} > 0 $ such
that $ A_n - z $ is invertible for $ n \geq n_{z_0,h} $ and $
|z-z_0| < \delta_{z_0,h} $. By compactness, $ U_0 $ can be
covered by finitely many balls of the form $ \{ |z-z_j| <
\delta_{z_j,h} \} $ and thus $ A_n - z $ is invertible for all $ z
\in U_0 $ and $ n \geq n_{h,U_0} := \max_{j} n_{h,z_j} $. The
balls can be chosen such that
$$ \sup_{n \geq n_{h,z_j}} \sup_{|z-z_j| < \delta_{h,z_j}}|| (A_n-z)^{-1} ||_{H^{s,\sigma}_{\rm sc}
\rightarrow H^{s+2,\sigma}_{\rm sc}} < + \infty $$ so the norms $
|| (A_n-z)^{-1} ||_{H^{s,\sigma}_{\rm sc} \rightarrow
H^{s+2,\sigma}_{\rm sc}}  $ are uniformly bounded with respect to
$ n \geq n_{h,U_0} $ and $ z \in U_0 $. Then
(\ref{convergenceabstraite}) follows  from the resolvent identity. \finpreuve

\bigskip

For $ k \geq 1 $ integer, to be fixed further on, we set
$$ f_z^k (\lambda) = (\lambda-z)^{- k} . $$
%In the following proposition, we shall justify the existence of such traces.
\begin{prop} \label{pourlimite} Let  $ U \subset \Ca $ an open subset and $ a = (a(h))_{h  \in (0,h_0] } $ be a family of $ S^{2,0} $
satisfying (\ref{elliptiquePDE}). Let $ b =(b(h))_{h \in (0,h_0]}
$ be a family of $ S^{m,\mu} $ with $ m \leq 2 $ and $ \mu < 0 $.
Assume that, for all $ h \in (0,h_0] $ and all $ z \in U $,
$$ A - z  : H^{2,0}_{\rm sc} \rightarrow L^2 (\Ra^d) $$
is invertible.

\noindent i) Let $ j \geq 1 $. Then, $ \frac{d^j}{d \varepsilon^j}
f^k_z(A+\varepsilon B)_{|\varepsilon = 0} $ is well defined and is
a linear combination of
\begin{eqnarray}
(A-z)^{-k_1}B (A-z)^{-k_2} \cdots B (A-z)^{-k_{j+1}}, \qquad k_1 +
\cdots + k_{j+1} = k+j \label{seulprouve}
\end{eqnarray}
with $ k_1 , \ldots ,k_{j+1} \geq 1 $. Furthermore, if
\begin{eqnarray}
   j (m-2) - 2 k < - d    \qquad \mbox{and} \qquad j \mu < - d ,
  \label{conditiontracej}
\end{eqnarray}
 each  operator of the form (\ref{seulprouve}) is of trace class in $ L^2 (\Ra^d) $.

\noindent ii)
 Assume in addition that, for all $ h \in (0,h_0] $ and all $
z \in U $,
$$   A + B - z : H^{2,0}_{\rm sc} \rightarrow L^2 (\Ra^d) $$
is invertible. Then
\begin{eqnarray}
 f_z^k (A + B) - f_z^k (A) - \sum_{j=1}^{p-1} \frac{1}{j!}
\frac{d^j}{d\varepsilon^j} f_z^k (A + \varepsilon B)_{|\varepsilon
= 0} \label{fausseexpression}
\end{eqnarray}
is well defined and is a linear combination of
\begin{eqnarray}
(A+B-z)^{-k_1} B (A-z)^{-k_2} \cdots B (A-z)^{-k_{p+1}}, \qquad
k_1 + \cdots + k_{p+1 } = k + p  \label{combinaisonlineaire}
\end{eqnarray}
with $ k_1 , \ldots , k_{p+1} \geq 1 $. If
\begin{eqnarray}
  p (m-2) - 2 k < - d    \qquad \mbox{and} \qquad p \mu < - d \label{conditiontracep}
\end{eqnarray} then each operator of the form (\ref{combinaisonlineaire}) is
trace class on $ L^2 (\Ra^d) $.
\end{prop}

First recall that  from the standard estimate
$$ || \scal{x}^{ -s} \scal{h D}^{-\sigma} ||_{\rm tr} \leq C h^{-d} , \qquad h \in (0,h_0] , $$
we have:
\begin{lemm} \label{lemmeclassiquetrace} For all $ s > d $ and $ \sigma > d $, the injection $
H^{s,\sigma}_{\rm sc} \hookrightarrow L^2 (\Ra^d) $ is trace class
with norm $ {\mathcal O}(h^{-d}) $.
\end{lemm}

%\noindent  {\it Proof.} This follows from the standard estimate
%$$ || \scal{x}^{ -s} \scal{h D}^{-\sigma} ||_{\rm tr}= {\mathcal O}(h^{-d}) , \qquad h \in (0,h_0] . $$
%\finpreuve

\bigskip

\noindent {\it Proof of Proposition \ref{pourlimite}.} That $
\frac{d^j}{d \varepsilon^j} f^k_z(A+\varepsilon B)_{|\varepsilon =
0} $ is well defined follows directly from Proposition
\ref{lemmeprop} {\it i)}, as well as its expression for $ k = 1$
which is given by (\ref{deriveeTaylor}). The formula for $ k \geq
2 $ is obtained by applying $
\partial_z^{k-1} $ to (\ref{deriveeTaylor}), using 
\begin{eqnarray}
 (k - 1)! (\lambda-z)^{-k} = \partial_z^{k-1} (\lambda-z)^{-1} ,
 \label{pourderiverresolvente}
\end{eqnarray}
and the smoothness of (\ref{maplocal}).
  By Proposition \ref{inverse},  each operator of the form
(\ref{seulprouve}) is bounded from $ L^2 (\Ra^d) $ to $
H^{j(2-m)+2k,-j\mu}_{\rm sc} $ thus is trace class  by
(\ref{conditiontracej}) and Lemma \ref{lemmeclassiquetrace}. This
completes the proof of {\it i)}. The proof of {\it ii)} is
completely similar once observed that, for $ k = 1 $,
(\ref{fausseexpression}) equals
$$ (-1)^p (A+B-z)^{-1} \left( B (A-z)^{-1} \right)^{p} \ ,$$
which is obtained
% by induction on $ p $
 using
(\ref{deriveeTaylor}). \finpreuve

\bigskip

\noindent {\bf Conclusion.} Under the assumptions of Proposition
\ref{pourlimite} {\it ii)},   the following expression is well
defined:
\begin{eqnarray}
T_p^k (A, B , z) := \mbox{tr} \left( f^k_z (A + B) - f^k_z (A) -
\sum_{j=1}^{p-1} \frac{1}{j!} \frac{d^j}{d\varepsilon^j} f^k_z (A
+ \varepsilon B)_{|\varepsilon = 0}  \right) , \label{Taylorkp}
\end{eqnarray}
(with the usual convention that $ \sum_{j=1}^{p-1} \equiv 0 $ if $
p = 1$) provided that (\ref{conditiontracep}) holds, thus in particular for 
$$ k > d/2 \qquad \mbox{and} \qquad p \mu < - d . $$
If in addition $ (a(h))_{h \in (0,h_0]} \in {\mathcal B} $ as in Proposition \ref{inverse}, we  have the following bound,
\begin{eqnarray}
\left| T_p^k (A,B,z) \right| \leq C h^{-d} \left( 1 + ||
(A-z)^{-1} ||_{L^2 \rightarrow H^{2,0}_{\rm sc}} + || ( A + B - z
)^{-1} ||_{L^2 \rightarrow H^{2,0}_{\rm sc}} \right)^N ,
\label{bornefinale}
\end{eqnarray}
for some $ C, N > 0  $ independent of $ h \in (0,h_0] $ and $ z
\in U $, using (\ref{normrestriction}), (\ref{combinaisonlineaire}) and Lemma
\ref{lemmeclassiquetrace}.

\section{Resonances} \label{OpRe}
\setcounter{equation}{0}
\subsection{The analytic distortion method}
In this subsection, we recall the definition of  resonances by the
analytic distortion method after Sj\"ostrand-Zworski. We also
collect additional results that will be necessary for our
applications.

We first recall the definition of    a maximal totally real
manifold  $ \Gamma \subset \Ca^d  $ parametrized by $ \kappa :
\Ra^d \rightarrow \Ca^d  $. By this it is meant
 that $ \kappa : \Ra^d \rightarrow \kappa (\Ra^d) = \Gamma  $
is a diffeomorphism (between real manifolds) such that
\begin{eqnarray}
T_{\kappa(x)} \Gamma \cap i (  T_{\kappa(x)} \Gamma   )  = \{ 0\} ,
\qquad x \in \Ra^d . \nonumber
\end{eqnarray}
Equivalently this means that, for all $x$, $  (\partial_1
\kappa(x) , \ldots ,
\partial_d \kappa (x) , i \partial_1 \kappa (x) , \ldots , i \partial_d
\kappa(x))  $ is a basis of  $ \Ca^d  $ viewed as a real vector
space, or that $ (\partial_1 \kappa (x) , \ldots , \partial_d \kappa(x))
$ is a basis of $ \Ca^d $ as a complex vector space, so the fact
that $ \Gamma $ is totally real simply means that
\begin{eqnarray}
  \mbox{det} \left( \frac{\partial \kappa(x)}{\partial x}
\right) \ne 0 , \qquad x \in \Ra^d . \label{totallyreal}
\end{eqnarray}
%Given such a $ \kappa  $, we consider an almost analytic extension $
%\tilde{\kappa} : \Ca^n \rightarrow \Ca^n  $ of $ \kappa $, ie a smooth map
%such that  $ \overline{\partial} \tilde{\kappa}  $ vanishes to infinite order on $
%\Ra^n  $  and $ \tilde{\kappa}_{|\Ra^n} =\kappa $.
Then, to any differential operator
$$ P =
 \sum_{|\alpha| \leq m} a_{\alpha}(x) D^{\alpha} , $$ with
 coefficients that are smooth on $ \Ra^d  $ and
holomorphic in some neighborhood of $  \Gamma \cap \left( \Ca^d
   \setminus \Ra^d  \right)  $ (typically a sector of the form $ \Sigma (\theta_0,R_0, \epsilon_0) $), we
 can associate the operator
\begin{eqnarray}
 {\mathcal A}_{\kappa} P := \sum_{|\alpha| \leq m} a_{\alpha}
(\kappa(x)) \left( (^t\partial_x \kappa (x) )^{-1} D \right)^{\alpha}
. \label{parametrisation}
\end{eqnarray}
%Of course, such a definition is independent of the choice of the almost
%analytic extension $ \tilde{\kappa}  $.

 The analytic
distortion method is as follows. Given $ R_1 > 0 $ and $
\epsilon_1  > 0  $, we can find a non decreasing smooth function $
\phi : \Ra^+ \rightarrow \Ra $ such that
\begin{eqnarray}
 \phi (t) & = & 0  \qquad t \leq R_1, \label{Rvariete} \\
  \phi (t) & = & 1  \qquad t \gg 1 , \label{grandt} \\
%  1 + i t \theta \varphi^{\prime}(t) & \ne & 0, \qquad t > 0 , \ \theta \in [0, \pi / 2 ] , \nonumber \\
0 \leq \  t \theta \phi^{\prime}(t)  & \leq & \epsilon_1 , \qquad
t > 0, \ \theta \in [0,\pi] , \label{epsilonvariete1}
\end{eqnarray}
and the latter condition implies, by possibly considering $ \phi $
associated with a smaller $ \epsilon_1 $, that we can additionally
assume
\begin{eqnarray}
0 \leq \mbox{arg} ( 1 + i t \theta \phi^{\prime}(t) )  \leq
\epsilon_1 , \qquad t > 0, \ \theta \in [0,\pi] .
\label{epsilonvariete2}
\end{eqnarray}
We assume in the sequel that, for each $ \epsilon_1  > 0$ (small
enough) and $ R_1 > 0 $ (large enough), a function $ \phi $
satisfying (\ref{Rvariete}), (\ref{grandt}),
(\ref{epsilonvariete1}) and (\ref{epsilonvariete2}) has been
chosen. Then the function
$$ f_{\theta}(t) = e^{ i \phi (t)
\theta} t , \qquad t \in \Ra^+ , $$ satisfies
\begin{eqnarray}
f_{\theta}(t) = t \ \ \mbox{for} \ \ t \leq R_1, \qquad
f_{\theta}(t) = e^{i \theta} t \ \ \mbox{for} \ \ t \gg 1 , \qquad
\partial_t f_{\theta} \ne 0 \label{ellipticstable} \\
0 \leq \mbox{arg} (f_{\theta}(t)) \leq \theta , \qquad
\mbox{arg}(f_{\theta}(t)) \leq \mbox{arg}(\partial_t
f_{\theta}(t)) \leq \mbox{arg}(f_{\theta}(t)) + \epsilon_1 .
\label{utilite}
\end{eqnarray}
Using this function, we can now define $ \kappa_{\theta} : \Ra^d
\rightarrow \Ca^d $ and $ \Gamma_{\theta} $ by
\begin{eqnarray}
 \kappa_{\theta} (x) = f_{\theta}(|x|) \frac{x}{|x|} = e^{i \theta \phi (|x|)}x , \qquad
\Gamma_{\theta} =  \kappa_{\theta}(\Ra^d) . \label{diffetvariete}
\end{eqnarray}
Notice that,
\begin{eqnarray}
\partial_x \kappa_{\theta} (x) = e^{i \theta \phi (|x|)} \left(
I_d  + i \theta |x| \phi^{\prime}(|x|) \frac{x \otimes x}{|x|^2}
\right), \label{secteur}
\end{eqnarray}
 thus (\ref{totallyreal}) holds, at least for $ \epsilon_1 $ small enough.  Now, if $ P $ is a differential operator whose coefficients
can be continued analytically to $ \Sigma
(\theta_0,R_0,\epsilon_0) $, by choosing $  \epsilon_1   $ small enough and
$$ R_1 > R_0 , \qquad  0   \leq \theta \leq \theta_0 , $$
 we can define the following
differential operator on $ \Ra^d $
\begin{eqnarray}
 P (\theta ) := {\mathcal A}_{\kappa_{\theta}} P , \label{notationdistortion}
\end{eqnarray}
 with  $ {\mathcal
A}_{\kappa_{\theta}} $ defined by (\ref{parametrisation}) and
(\ref{diffetvariete}).

\medskip

\noindent {\bf Remark.} %We use the simplified notation $ P
%(\theta) $ for $  {\mathcal A}_{\kappa_{\theta}} P $, however 
The
reader should keep in mind that operators of the form $ P (\theta) $ depend not only on
$ \theta $ (and $h \in (0,h_0]$ below) but also on the parameters
$ R_1 $ and $ \epsilon_1 $ (and also on the choice of the function
$ \phi $), although this dependence is omitted in the notation.

\medskip

\begin{defi} Let $V \in {\mathcal V}_{\rho} (\theta_0, R_0, \epsilon_0) $.
The pair $ ( R_1 , \epsilon_1 ) \in \Ra_+^2 $ is said to be
\underline{Fredholm admissible} for $H_0  + V $  if, for all $
\theta \in [0,\theta_0 ] $, the following hold:

\noindent i) for all  $ h \ll 1$ and all $ z \in \Ca \setminus e^{-2i\theta}[0,+\infty)
$,
$$ H_0 (\theta) + V (\theta) - z : H^2(\Ra^d) \rightarrow L^2 (\Ra^d) \ \ \mbox{is a Fredholm operator of index} \ 0, $$
ii) the principal symbol, in the classical sense, $
p_{\theta}^{\rm cl} $ of $ H_0 (\theta) + V (\theta) $ is
elliptic, ie for some $ C \geq 1 $
$$ |p_{\theta}^{\rm cl}(x,\xi)| \geq C^{-1}|\xi|^2 , \qquad (x,\xi) \in \Ra^{2d} . $$

   Here $ H_0 (\theta) $ and
$ V (\theta) $ are defined by (\ref{notationdistortion}) with $
\kappa_{\theta} $ given by (\ref{diffetvariete}).
% and such that $
%\mbox{det}(\partial_x \kappa_{\theta}) \ne 0 $ on $ \Ra^d $.
\end{defi}

\begin{prop} \label{Fredholmprop} Let $  (V_{\iota})_{\iota \in I}
$ be  bounded family of  $ {\mathcal V}_{\rho}
(\theta_0,R_0,\epsilon_0) $. We can find $ \overline{R}_1
> 0 $, $ \overline{\epsilon}_1 > 0 $ and $ \overline{C} > 0 $ such that, for all $ \iota \in I $,
any  $( R_1, \epsilon_1 ) \in [\overline{R}_1 , + \infty) \times
(0, \overline{\epsilon}_1] $ is Fredholm admissible for $ H_0 +
V_{\iota} $, with  constant $\overline{C} $ in ii). More
explicitly
\begin{eqnarray}
| p_{\iota,\theta}^{\rm cl}(x,\xi) | \geq \overline{C}^{-1}
|\xi|^2 , \label{controleellipticite}
\end{eqnarray}
uniformly with respect to $ \epsilon_1 \in
(0,\overline{\epsilon_1} ]$, $ R_1 \geq \overline{R_1}$, $ \theta
\in [0,\theta_0] $ and $ \iota \in I $. In addition, we may also
assume that, for all $ \theta \in [0,\theta_0] $,
\begin{eqnarray}
- 2 \theta - 3 \epsilon_1 \leq \emph{arg} \left(
p_{\iota,\theta}^{\rm cl}(x,\xi) \right) \leq \epsilon_1 , \qquad
x \in \Ra^d, \ \xi \in \Ra^d \setminus 0 ,
\label{controleargument}
\end{eqnarray}
uniformly with respect to $ \iota \in I $.
\end{prop}

 Proposition \ref{Fredholmprop} is proved, for
a single $ V $,  in the lecture notes \cite[Lemma 7.3]{SjBook} in
the more general framework of {\it black box perturbations}. Its
extension to a bounded family of $ {\mathcal V}_{\rho}
(\theta_0,R_0,\epsilon_0) $ involves no new argument and we
therefore omit the proof. The reason for considering a bounded
family in $ {\mathcal V}_{\rho} (\theta_0,R_0,\epsilon_0) $ is
that we shall approximate $ V \in {\mathcal V}_{\rho}
(\theta_0,R_0,\epsilon_0) $ by a sequence $ V_n \in {\mathcal
V}_{\bar{d}} (\theta_0,R_0,\epsilon_0) $, with $ \bar{d} > d $,
and use a certain deformation along $ \kappa_{\theta}(\Ra^d) $. It
will be important that $ \kappa_{\theta} $ (which depends on $
\epsilon_1 $ and $ R_1 $) can be chosen independently of $n$.

The Fredholm admissibility is important to define the resonances
as we shall see below. In the case of a single $ V$, the first
part of Proposition \ref{Fredholmprop} simply states that this
condition is fulfilled for  $ H_0 + V $. The additional uniform
estimates (\ref{controleellipticite}) and (\ref{controleargument})
will be useful later on to prove some resolvent estimates.

\bigskip

The definition of resonances relies on the following theorem.

\begin{theo} \label{admis0} (\cite{SjostrandZworski,Sjostrace,SjBook}) Let $ 0  < \theta_0 < \pi  $ and $
 V \in {\mathcal V}_{\rho} (\theta_0,R_0,\epsilon_0) $.
 Assume that we are given
$R_1 > 0 $
%in $ \reff{ellipticstable} $
and $ \epsilon_1 > 0 $ which are Fredholm admissible. Then, for
all $ h \ll 1 $ and all $ z \in \Omega$, we have

\noindent (i) $ z \in \sigma ( H_0 (\theta ) + V (\theta)) $ if
and only if $ \emph{ker}(H_0 (\theta) + V (\theta)-z) \ne 0 $.

\noindent (ii) For all $ 0 \leq \theta_1 \leq \theta_2 \leq
\theta_0 $, if $ z \in \Ca \setminus e^{-2 i [ \theta_1 , \theta_2
]}[0,+\infty) $ then
$$ \emph{dim ker} (H_0 (\theta_1) + V (\theta_1) - z) = \emph{dim ker}
(H_0 (\theta_2) + V (\theta_2) - z) . $$
\end{theo}

The first statement is an immediate consequence of the fact that
the operator has a zero index. The second one requires a non
trivial analytic deformation result, which uses the analyticity of
the coefficients of $V$ near infinity.

 Let us recall the main consequence of Theorem \ref{admis0}.

 First, if $ 0 \leq \theta \leq
\theta_0 < \pi $ and  $ 0 < \epsilon < 2 \pi - 2 \theta_0 $, then
for   all $ h \ll 1 $ and all $ z \in e^{i(0,\epsilon)}
(0,+\infty) $,
\begin{eqnarray}
 H_0 (\theta) + V (\theta) - z : H^2 (\Ra^d) \rightarrow L^2 (\Ra^d) \ \mbox{is an
isomorphism}  . \label{sansresonance}
\end{eqnarray}
Furthermore, by analytic Fredholm theory, one can show that the
spectrum of $ H_0 (\theta) + V (\theta) $ is discrete in $ \Ca
\setminus e^{- 2 i  \theta} [ 0 , + \infty ) $. The part {\it
(ii)} guarantees that, if $ \theta^{\prime} > \theta $, the
eigenvalues of  $ H_0 (\theta) + V (\theta) $ and
 $ H_0 (\theta^{\prime}) +
V (\theta^{\prime}) $ coincide on $ e^{-2i [0,\theta)}(0,+\infty) $ and this makes the following definition natural.
\begin{defi} Given $ \Omega $ satisfying (\ref{localiseOmega}),
the set of resonances of $ H_0 + V $ in $\Omega$ is
% (and (\ref{intervalle}),
$$ \emph{Res} (H_0 + V , \Omega) = \Omega \cap  \sigma (H_0 (\theta_0) + V (\theta_0))
\cap e^{-i[0, 2 \theta_0)}(0,+\infty) . $$
\end{defi}
Recall that $ \mbox{Res} (H_0 + V , \Omega) $ is finite (for each
$ h $).

 By analytic Fredholm theory again, for
any  $ w \in \mbox{Res}(H_0+V,\Omega)$, the operator
\begin{equation}\label{defprojecteur}
 \Pi_{\theta,w} = \frac{i}{2  \pi} \int_{\gamma(w)} (H_0 (\theta) + V (\theta) - z )^{-1}dz
 \end{equation}
 is of finite rank, if  $ \gamma (w) $ a small enough contour enclosing
 $ w $ and this allows to state the following definition.
\begin{defi} \label{defimultiplicite}
The  multiplicity of  $ w $ is the rank of $ \Pi_{\theta,w} $.
\end{defi}
This definition is independent of $ \theta $ in the sense that we
get the same rank if $ \theta $ is replaced by some larger $
\theta^{\prime} $ (smaller than $  \theta_{0} $).

We conclude this subsection with the following elementary resolvent
estimates.

\begin{prop} \label{Omegadeltaresolvente} Let $ \Omega $  be  satisfying (\ref{localiseOmega}) and let $ \Omega^+_{\delta} :=
\Omega^+ \cap \{ \emph{Im}(z) \geq \delta \} $ (see (\ref{Omegasuperieur})) with $ \delta $
small enough to be non empty. Let $ (V_{\iota})_{\iota \in
I} $ be a bounded family of $ {\mathcal V}_{\rho}
(\theta_0,R_0,\epsilon_0) $. Then,  for all $ \epsilon_1
> 0 $ small enough, we can choose  $ R_1 > 0 $ as large as we want such that
\begin{eqnarray}
|| (H_0(\theta_0)+V_{\iota}(\theta_0)-z)^{-1}||_{L^2 \rightarrow
H^{2,0}_{\rm sc}} \lesssim 1, \qquad h \ll 1 , \ z \in
\Omega^+_{\delta}, \ \iota \in I \ . \label{borneresolventedelta}
\end{eqnarray}
\end{prop}

\noindent {\it Proof.} Denote by $ p_{\iota}(x,\xi,h) $ the full
Weyl symbol of $ H_0 + V_{\iota} $, which is then real on $
\Ra^{2d} $ and of the form
$$ p_{\iota}(x,\xi,h) = p^{\rm cl}_{\iota}(x,\xi) +
a_{\iota} (x,\xi,h) , $$ with $ a_{\iota} (.,h)$ polynomial of
degree $ \leq 1 $ in $ \xi $ with coefficients bounded in $
{\mathcal C}_{\rho}(\theta_0,R_0,\epsilon_0) $. Setting
$$ p_{\iota,\theta_0}(x,\xi,h) = p_{\iota} \Big(\kappa_{\theta_0}(x), (^t \partial_x \kappa_{\theta_0}(x))^{-1} \xi , h \Big), $$
we then have
$$  H_0 (\theta_0) + V_{\iota} (\theta_0) = O \! p^w_h (p_{\iota,\theta_0}) + h O \! p_h^w (b_{\iota,\theta_0}(h)) $$
for some symbol $ b_{\iota,\theta_0}(h) $ which, for fixed $
\epsilon_1 $ and $ R_1 $, is bounded in $ S^{1,0} $ as $ h $ and $
\iota $ vary. We thus only need to show that, for  $ \epsilon_1
> 0 $ small enough and  $ R_1 > 0 $ large enough,
\begin{eqnarray}
|p_{\iota,\theta_0}(x,\xi,h)-z| \gtrsim 1 , \qquad h \ll 1 , \ z
\in \Omega^+_{\delta}, \ \iota \in I . \label{minoration1}
\end{eqnarray}
The result then follows from the standard construction of a
semiclassical parametrix, yielding the  invertibility of $
H_0(\theta_0)+V_{\iota}(\theta_0)-z $ for $ h \ $ small enough
(uniformly with respect to $z$ and $\iota$) as well as the bound
(\ref{borneresolventedelta}). Let us prove (\ref{minoration1}).
Using (\ref{controleellipticite}), we can choose $ C_0 > 0 $ large
enough, independent of $ 0 < \epsilon_1 \leq
\overline{\epsilon}_1$, $ R_1 \geq \overline{R}_1 $, $ x \in \Ra^d
$, $ h \ll 1 $ and $ \iota \in I $ such that
$$ | p_{\iota,\theta_0}(x,\xi,h) | \geq 1 + \max_{\overline{\Omega}}|z|,
\qquad  |\xi| \geq C_0 ,$$ since $ | (p_{\iota,\theta_0} -
p_{\iota,\theta_0}^{\rm cl})(x,\xi,h) | \lesssim \scal{\xi} $, uniformly with respect to $ h, \iota,\epsilon_1,R_1 $.
 Using (\ref{controleargument}), if $ \epsilon_1 > 0 $ and $  \delta^{\prime} > 0 $ are  small enough, we also have
$$ | p_{\iota,\theta_0}^{\rm cl}(x,\xi) - z | \geq \delta^{\prime}, \qquad x,\xi \in \Ra^d , \ z \in \Omega^+_{\delta}.   $$
Then, once such  $ \epsilon_1 $ and $ \delta^{\prime} $ have been chosen, we have, for all $ R_1 $ large enough,
$$  | a_{\iota} \Big(\kappa_{\theta_0}(x), (^t \partial_x \kappa_{\theta_0}(x))^{-1} \xi , h \Big) | \leq
 \frac{\delta^{\prime}}{2} ,
\qquad |x| \geq R_1 , \qquad |\xi| \leq C_0 ,
 $$
 since the coefficients of $ a_{\iota} $ decay like $ \scal{x}^{-\rho}
 $ in $ \Sigma (\theta_0,R_0,\epsilon_0) $ uniformly with respect
 to $ h $ and $\iota $.
% Recall that $ a $ is real on $ \Ra^{2d} $ so we also have
%$$ a (\kappa_{\theta}(x), ^t \partial_x \kappa_{\theta}(x)^{-1} \xi , h ) \in \Ra $$
It is then straightforward to check that
%$$  |p_{\iota,\theta_0}(x,\xi,h)-z| \geq \min \left( 1 , \delta/2 \right),  $$
%ie 
(\ref{minoration1}) holds since  $ p_{\iota,\theta_0} $
is real for $ |x| \leq R_1 $. \finpreuve

\bigskip

In the next proposition, we prove an exponential bound for the resolvent of $ H_0 (\theta) $. The latter can be used with Theorem \ref{theoremeavecborne} to obtain an exponential upper bound on $ \partial_z \varphi_p (z,h) $, when $ p \geq 3 $. Let us recall that, since $ H_0 = - h^2 \Delta $ has no resonances away from $0$, $ (H_0 (\theta)-z)^{-1} $ is well defined for all $ z \in \Omega $ (see \cite{SjostrandZworski}).

For simplicity, we only consider the case where $ \theta_0 < \pi /2 $ and $ d \geq 3 $.

\begin{prop} \label{propestires} Assume that $ \theta_0 < \pi /2 $ and that $ d \geq 3 $.
 Let $ \Omega $ be a simply connected open set satisfying (\ref{localiseOmega}). Then, if
 $ \epsilon $ (in (\ref{localiseOmega})) and $ \epsilon_1 $ (in (\ref{epsilonvariete1}))
 are small enough, we have
\begin{eqnarray}
|| (H_0(\theta_0)-z)^{-1}||_{L^2 \rightarrow
H^{2,0}_{\rm sc}} \lesssim e^{Ch^{-1}}, \qquad h \ll 1 , \ z \in
\Omega \ . \label{borneresolvente}
\end{eqnarray}
\end{prop}

%Notice that choosing $ \theta_0 < \pi / 2 $ is probably not optimal. Our point with this proposition is simply to give an example of 
%relatively explicit upper bound for $ (H_0(\theta)-z)^{-1} $.

\bigskip

\noindent {\it Proof.} %According to Proposition
%\ref{Omegadeltaresolvente} and to the simple connectedness of
%$\Omega$ it is sufficient to prove the estimate for
%$(H_0(\theta_0)-z)^{-k}$, $k >d/2+1$.
By (\ref{diffetvariete}) and (\ref{secteur}), the coefficients of $ H_0 (\theta) $ are
holomorphic with respect to $ \theta $ in a small neighborhood of
$ [0,\theta_0 ] $ and thus so is
\begin{eqnarray}
 \theta \mapsto \left( v,(H_0(\theta)-z)^{-1} u  \right),
 \label{weakaexpliciter}
\end{eqnarray}
for $ \theta $ in  a complex neighborhood of $ [0,\theta_0] $ and for all $ u , v \in C_0^{\infty}(\Ra^d) $, $ z \in \Omega $ and $
h \in (0,1] $. On the other hand, for $ i \theta \in \Ra $ small,
$$H_0(\theta)= U_{\theta} \, H_0 \, U_{\theta}^{-1},$$
with $U_{\theta} : L^2(\Ra^d) \rightarrow L^2(\R^d)$ the
isomorphism defined by $U_{\theta}(u)(x)= u(\kappa_{\theta}(x))$.
Since $ U_{\theta} $ maps $ H^2 (\Ra^d) $ into itself, we then have
$$ (H_0(\theta)-z)^{-1} =  U_{\theta} (H_0 -z)^{-1} U_{\theta}^{-1} , \qquad z \in \Omega^+ ,$$
and if we denote by $ {\mathcal R}(x-y,z,h) $ the Schwartz kernel
of $ (H_0 -z)^{-1} $
 we can rewrite (\ref{weakaexpliciter}) as
\begin{eqnarray}
  \int_{\Ra^{2d}} {\mathcal R}(
\kappa_{\theta}(x)-\kappa_{\theta}(y),z,h)
  u (y) \overline{v(x)} \mbox{det}(\kappa_{\theta}(y)) \ d x d y
  \label{doubleprolongement}
\end{eqnarray}
for $ i \theta \in \Ra $ small and $ z \in \Omega^+ $. Let us
recall that, for $ \mbox{Im}(z^{1/2}) > 0 $,
$$ {\mathcal R}( x-y,z,h) = \frac{i}{4 h^2 } \left( \frac{z^{1/2}}{2 \pi h |x-y|}
\right)^{\frac{d}{2}-1} H^1_{\frac{d}{2}-1}( z^{1/2}|x-y|/h) , $$
where the Hankel function $ H^1_{\nu}(Z) $ (with $ \nu =
\frac{d}{2}-1 $) is given by
$$ H_{\nu}^1 (Z) = \left( \frac{2}{\pi Z} \right)^{1/2} \frac{e^{i (Z - \frac{\nu}{2}\pi - \frac{\pi}{4})}}{
\Gamma (\nu + \frac{1}{2})} \int_0^{+\infty} e^{-s} \left(s (1 +
i s Z^{-1}/2) \right)^{\nu - \frac{1}{2}} ds ,
$$
using everywhere the determination of the square root defined on $
\Ca \setminus (-\infty,0 ] $ taking its values in $
e^{i(-\pi/2,\pi/2)} (0,+\infty) $ (see for instance section VII.7.2 of \cite{Watson}). The function $ H_{\nu}^1 $ is
holomorphic for $ Z \in e^{i(-\pi/2,\pi/2)} (0,+\infty) $, with
the following rough bound, for all $ 0 < \delta < \pi/2$,
\begin{eqnarray}
| H_{\nu}^1 (Z) |  \leq  C_{\delta}|Z|^{-1/2} e^{|\rm{Im}(Z)|}
\max \left(1 , |Z|^{\frac{1}{2}- \nu} \right) , \qquad
\mbox{arg}(Z) \in (\delta -\pi/2, \pi/2- \delta) .
\label{borneHankel}
% \\
%| H_{\nu}^1 (Z) | & \lesssim & \qquad |Z| \rightarrow \infty , \
%\mbox{arg}(Z) \in (-\pi/2,\pi/2)
\end{eqnarray}
Independently, by writing $ \varphi (x) = \phi (|x|) $, we have
$$ \kappa_{\theta}(x) - \kappa_{\theta}(y) = (x-y) \int_0^1 e^{i \theta \varphi(y + t (x-y))}
\Big(  i\theta \nabla \varphi (y + t (x-y)) \otimes (y + t (x-y)) + 1
\Big) dt , $$
 where $ | \theta \nabla \varphi (X) \otimes X| \lesssim
\epsilon_1 $ by (\ref{epsilonvariete1}) and $0 \leq \varphi(X)
\leq 1 $. Therefore, if $ \epsilon_1 $ and  $ \epsilon $ are small enough, there
exists $ \delta > 0 $ small enough such that
$$ z^{1/2}|\kappa_{\theta}(x)-\kappa_{\theta}(y)| := \left( z \scal{\kappa_{\theta}(x)-\kappa_{\theta}(y) ,
\kappa_{\theta}(x)-\kappa_{\theta}(y) } \right)^{1/2} \in
e^{i(\delta-\pi/2,\pi/2-\delta)}(0,+\infty) ,
$$
 for $ x \ne y  $, $ x,y \in \Ra^d $, $ z \in \Omega $ and $
\theta $ in a neighborhood of $ [0,\theta_0] $. Furthermore, the modulus of $ |\kappa_{\theta}(x)-\kappa_{\theta}(y)| / |x-y| $ is bounded from above and from below.
 This allows to continue
(\ref{doubleprolongement}) analytically with respect to $ \theta
\in [0,\theta_0] $ and then with respect to $ z \in {\Omega} $.
Using (\ref{borneHankel}) and the Schur Lemma, we deduce that, for
any $ \chi \in C_0^{\infty}(\Ra^d) $,
$$ || \chi ( H_0(\theta_0 )-z)^{-1} \chi ||_{L^2 \rightarrow
L^2} \lesssim e^{C h^{-1}} , \qquad z \in \Omega. $$ This easily
implies a similar $ L^2 \rightarrow L^2 $ bound on the whole
resolvent using the elementary estimate
$$
|| (e^{-2i\theta_0 } H_0-z)^{-1}||_{L^2 \rightarrow
H^{2,0}_{\rm sc}} \lesssim 1, \qquad z \in \Omega,
$$
 and two
applications of the resolvent identity yielding
\begin{eqnarray*}
  ( H_0(\theta_0 )-z)^{-1} & = &  (e^{-2i\theta_0 } H_0-z)^{-1} -
(e^{-2i\theta_0 } H_0-z)^{-1} V_0 (e^{-2i\theta_0 } H_0-z)^{-1} \\
& & + \ (e^{-2i\theta_0 } H_0-z)^{-1} V_0 ( H_0(\theta_0 )-z)^{-1}
V_0 (e^{-2i\theta_0 } H_0-z)^{-1} ,
\end{eqnarray*}
where $ V_0 := H_{0}(\theta_0) - e^{-2 i \theta_0}H_0 $ is a
compactly supported differential operator of order $2 $. The $ L^2
\rightarrow H^{2,0}_{\rm sc} $ bound then follows from the $ L^2 \rightarrow L^2 $ one by the resolvent identity between $ z_0 \in \Omega_{\delta}^+ $ and $z$, using (\ref{borneresolventedelta}). \finpreuve

\subsection{A deformation result}
We recall first the following result.
\begin{prop}[Sj\"ostrand \cite{Sjosequi}] \label{deformationSjtrace} Let $ \overline{d} > d $
and $ V \in {\mathcal V}_{\overline{d}} (
\theta_0,R_0,\epsilon_0 ) $. Let $ R_1 > 0
 $ and $ \epsilon_1 > 0 $ be Fredholm admissible for $ H_0 $ and $
 H_0 + V $. Then, if $ k > d/2 + 1 $,
 $$ \emph{tr} \left( (H_0 + V -z)^{-k} - (H_0 -z)^{-k} \right) =
 \emph{tr} \left( (H_0 (\theta) + V (\theta) - z)^{-k} - (H_0 (\theta) -z)^{-k} \right) , $$
 for all $ \theta \in [ 0 , \theta_0 ] $ and all $ z \in \Omega^+
 $.
\end{prop}

In the next proposition, we simply state that the above invariance of the trace by analytic distortion
still holds for the regularized traces of the form
(\ref{Taylorkp}).

\begin{prop} \label{deformationSjp} Let $ p \in \Na $ and $ \rho > 0
$ such that $ \rho > d/p$. Let $ V \in {\mathcal V}_{\rho} (\theta_0,R_0,\epsilon_0) $.
Then, if $ \epsilon_1 $ is small enough, $ R_1 $ is large enough
and $ k > d/2 + 1 $, we have
$$ T_p^k \left( H_0 , V , z \right) = T_p^k \left( H_0 (\theta), V (\theta) ,z \right) , $$
 for all $ \theta \in [ 0 , \theta_0 ] $ and all $ z \in \Omega^+ $.
\end{prop}

As the reader may guess, this proposition is a fairly elementary consequence of Proposition
\ref{deformationSjtrace}, approximating $ V $ by a sequence $ V_n
\in {\mathcal V}_{\bar{d}}
(\theta_0,R_0^{\prime},\epsilon_0^{\prime} ) $ with $ \bar{d}>d $.

%Since we shall need such a
%sequence at several places in this paper, we state the following
%lemma.

\begin{lemm} \label{approxdistordu} Let  $ V \in {\mathcal V}_{\rho}
(\theta_0,R_0,\epsilon_0) $. Let $ \bar{d}>d $. We can find $
R^{\prime}_0
> R_0 $, $ 0 < \epsilon_0^{\prime} \leq \epsilon_0 $ and a sequence $ (V_n)_{n \geq 1} \in {\mathcal
V}_{\bar{d}} (\theta_0,R_0^{\prime},\epsilon_0^{\prime}) $,
bounded in $ {\mathcal V}_{\rho}
(\theta_0,R_0^{\prime},\epsilon_0^{\prime}) $ such that, for all $
\rho^{\prime} < \rho $ and all $ s,\sigma \in \Ra $,
\begin{eqnarray}
 || V_n - V ||_{H_{\rm sc}^{s,\sigma} \rightarrow H_{\rm
sc}^{s-2,\sigma+\rho^{\prime}}} \rightarrow 0, \qquad n
\rightarrow \infty , \label{compacitelimite}
\end{eqnarray}
 for all $ h \ll 1 $.
%$$ \chi_n (x)  = \exp \left( - e^{- i \theta_0} |x|^2 / n \right),
%\qquad n \gg 1 , \ x \in \Ra^d .
%$$ If $ \epsilon_0 $ is small enough, then, for $ n_0  $ large enough,
%$$ V_n := \chi_n V \overline{\chi}_n $$
%belongs to $ {\mathcal V}_{\overline{d}} $ for all $ \overline{d}
%> d $,  the sequence $ (V_n)_{n \geq n_0} $ is bounded in $ {\mathcal V}_{\rho} (\theta_0,R_0,\epsilon_0) $ and
%
%In addition, if $ V $ is dilation analytic, so are the $ V_n $.
\end{lemm}

\noindent {\it Proof.} Choose first a determination of $ Z \mapsto
Z^{1/4} $ for $ Z \in \Ca \setminus e^{ 2 i \theta^{\prime}_0 }[ 0 , + \infty ) $, with $ \theta_0 < \theta^{\prime}_0 < \pi $. We may assume that it is positive on $ \Ra^+ $.  Choose also $ \chi \in
C_0^{\infty}(\Ra^d) $ such that $ 0 \leq \chi \leq 1 $, $ \chi (x)
\equiv 1 $ for $ |x| \leq R_0^{\prime}/2 $, and $ \chi (x) = 0 $
for $ |x| \geq R_0^{\prime} $. We then define
$$ \chi_n (x) = \chi (x) + (1-\chi (x)) \exp \left( -( x^2 )^{1/4} / n \right), \qquad n \geq 1 , $$
with $ x^2 = x_1^2 + \cdots +  x_d^2 $, and 
$$ V_n = \chi_n V \chi_n . $$
If $ R_0^{\prime} $ is large enough, the coefficients of $V_n $
are then such that (\ref{symboleprincipal}), (\ref{ellipticite})
and (\ref{autoadjoint}) hold, with $ c $ independent of $n$ in
(\ref{ellipticite}), and (\ref{compacitelimite}) is elementary.
Furthermore, if $ \epsilon_0^{\prime} $ is small enough $ x
\mapsto \exp \left(- ( x^2 )^{1/4} / n \right) $ has an analytic
continuation to $ \Sigma (\theta_0,R_0^{\prime},
\epsilon_0^{\prime} ) $ where it is uniformly bounded with respect
to $ n \geq 1 $. Therefore $ (V_n)_{n \geq 1} $ is bounded in $
{\mathcal V}_{\rho} (\theta_0,R_0^{\prime},\epsilon_0^{\prime}) $. Also,
it  clearly belongs to $ {\mathcal V}_{\bar{d}}
(\theta_0,R_0^{\prime},\epsilon_0^{\prime}) $ since, if $ x = t e^{i \theta} \omega $ with $ t \gg 1 $, $ \omega $ close to $ {\mathbb S}^{d-1} $
and $ \theta \in [0,\theta_0] $, we then have $ \mbox{Re} \left((x^2)^{1/4} \right) \gtrsim t^{1/2} \cos (\theta/2) \gtrsim t^{1/2} $. \finpreuve

\bigskip

\noindent {\it Proof of Proposition \ref{deformationSjp}.}  By
Proposition \ref{Fredholmprop},  for all $ R_1 $ large enough and
all $ \epsilon_1 $ small enough, $ (R_1,\epsilon_1) $ is Fredholm
admissible for  $ \varepsilon V_n $
 and  $ \varepsilon V $, for all $n \geq 1$ and $ \varepsilon \in [0,1]
 $. Using Proposition \ref{deformationSjtrace} with $ R^{\prime}_0 $ and $ \epsilon_0^{\prime} $, we then have
$$ \mbox{tr} \left( (H_0 + \varepsilon V_n -z)^{-k} - (H_0 -z)^{-k} \right) =
 \mbox{tr} \left( (H_0 (\theta) + \varepsilon V_n (\theta) - z)^{-k} - (H_0 (\theta) -z)^{-k}
 \right) $$
and the latter can be differentiated with respect to $ \varepsilon
$ using Proposition \ref{lemmeprop} since the operators inside the
trace are smooth with respect to $ \varepsilon $, in the trace
norm. This is easily seen, for instance for the left hand side,
by writing the operator inside the trace as a linear combination
of operators of the form
$$ (H_0 + \varepsilon V_n -z)^{-k_1} \varepsilon V_n (H_0 -z)^{-k_2}, \qquad k_1 + k_2 = k + 1 . $$
 Therefore,
$$ T_p^k \left( H_0 , V_n , z \right) = T_p^k \left( H_0 (\theta), V_n (\theta) ,z \right) $$
gives the result by letting $ n$ go to $ \infty $, using
(\ref{compacitelimite}) with $ \rho^{\prime} $ such that $ p
\rho^{\prime} > d $, Propositions \ref{lemmeprop} and
\ref{pourlimite}. \finpreuve

\subsection{The main tool of Sj\"ostrand's trace formula} \label{pourlafactorisation}

\begin{prop}
%[Sj\"ostrand]
 \label{sjos2}  Let $ \Omega $ be an open
subset satisfying
 (\ref{localiseOmega}) with $ 0 <  \theta_0  < \pi $ and $ 0 < \epsilon < 2 \pi - 2 \theta $. Let $ V \in {\mathcal
V}_{\rho}(\theta_0, R_0, \epsilon_0) $ with $\rho >0$.
 Then, we can fix $ h_1,
\epsilon_1 $ small enough and $ R_1 $ large enough such that
there exists a family of finite rank operators $ (K_{\varepsilon}
(\theta_0))_{0 < h \leq h_1,  \varepsilon \in [0,1]} $  with the
following properties:
\begin{eqnarray}
%||K_{\theta}||_1 +
 \emph{rank}(K_{\varepsilon} (\theta_0))
  \lesssim   h^{-d}  , \qquad \qquad \qquad \label{rank}
 \\
 || ( H_0 (\theta_0) + \varepsilon V (\theta_0)
+ K_{\varepsilon} (\theta_0) - z)^{-1}||_{L^2 \rightarrow
H^{2,0}_{\rm sc}}  \lesssim  1  , \label{pasdeparametrix}
\end{eqnarray}
for all $ h \in (0,h_1]$, $ z \in \Omega $, $ \varepsilon \in
[0,1]$. For all $ N, s, \sigma \in \Ra $ and $ k \in \Na $
\begin{eqnarray}
 || \partial_{\varepsilon}^k K_{\varepsilon} (\theta_0)   ||_{ H^{s,\sigma}_{\rm sc}\rightarrow
H^{N,N}_{\rm sc}}  \lesssim   1 , \qquad h \in (0,h_1] , \
\varepsilon \in [0,1] . \label{domaines}
\end{eqnarray}
 In addition,
 there exists $ \chi \in C_0^{\infty}(\Ra^d) $, independent of $h$ and $ \varepsilon $, such
that $ K_{\varepsilon } (\theta_0) = \chi K_{\varepsilon}
(\theta_0 ) \chi $.
\end{prop}

Note that %(\ref{rank}), 
(\ref{domaines}) and Lemma \ref{lemmeclassiquetrace}
%and (\ref{rank})
imply that
\begin{eqnarray}
|| \partial_{\varepsilon}^k K_{\varepsilon} (\theta_0)   ||_{ \rm
tr } & \lesssim &  h^{-d} , \qquad h \in (0,h_1] , \ \varepsilon
\in [0,1] . \label{traceslisses}
\end{eqnarray}

This proposition is essentially proved in \cite{Sjosequi, SjBook}. We
however recall the main argument of the proof to emphasize the
 dependence on $ \varepsilon $ which was not considered in those references.

\bigskip

\begin{lemm} \label{lemmeangulaire} For all $ \epsilon_1 > 0 $   such that
$ 2 \pi - 2 \theta_0 - 4 \epsilon_1 > \epsilon $ and $ \epsilon_1 < \epsilon $, and for all $ C \gg 1 $,
we can construct a smooth function $ F : D_F \rightarrow \Ca $,
with $ D_F $ a neighborhood of
$ e^{i[- 2 \theta_0 - 4 \epsilon_1 , \epsilon ]} [ 0 , + \infty
) $,
%\cup \{ z \ ; \ |z| \leq C^{-1}  \}, \label{domaineexplicite}
%\end{eqnarray}
such that
\begin{eqnarray}
 F (Z) = Z, \ \   \mbox{for} \ Z \ \mbox{such that} \ |Z| \notin [ C^{-1} , C ] \ \mbox{  or with argument close to} \
-2  \theta_0 , \label{donnesupport}
\end{eqnarray}
 and
\begin{eqnarray}
 | F (Z) - z | \gtrsim 1, \qquad  Z \in D_F , \ \ z \in \Omega .
 \label{donneinversible}
\end{eqnarray}
\end{lemm}

\noindent {\it Proof.} We can  define a function $ \mbox{arg}(Z) $
smooth  on $ e^{i(- 2 \theta_0 - 4 \epsilon_1^{\prime} ,
\epsilon^{\prime} )} ( 0 , + \infty ) $, with $ \epsilon_1^{\prime}
$ and $ \epsilon^{\prime} $ slightly larger that $ \epsilon_1 $
and $ \epsilon $ respectively, such that
$$ Z = |Z| \exp (i \mbox{arg}(Z) ) , \qquad \mbox{arg}(Z) \in (- 2 \theta_0 - 4
\epsilon_1^{\prime} , \epsilon^{\prime} ) . $$ Observe next that,
for some $ \theta < \theta_0 $ and $ r_2 > r_1 > 0 $,
\begin{eqnarray}
 \Omega \subset \{ z \in \Ca \ ; \ r_1 \leq |z| \leq r_2 , \ \ -
2 \theta \leq \mbox{arg}(z) \leq \epsilon  \} .
\label{ensemblededroite}
\end{eqnarray}
We next take $ C $ large enough so that  $ C^{-1} < r_1 < r_2 < C $ and choose $ \psi \in C_0^{\infty}(C^{-1},C) $ such that 
$ \psi \equiv 1 $
near $ [r_1, r_2] $.  For $\delta$ small enough, we also choose $
\Theta \in C^{\infty}(\Ra) $ non decreasing such that
$$ \Theta (\alpha) = \begin{cases} \mbox{const.} \geq - 2 \theta_0 - 2 \delta, & \mbox{if} \
\alpha < - 2 \theta_0 - 2 \delta \\
\alpha , & \mbox{if} \ |- 2  \theta_0 - \alpha | \leq
\delta \\
\mbox{const.} \leq - 2 \theta_0 + 2 \delta , & \mbox{if} \ \alpha
> - 2 \theta_0 + 2 \delta
\end{cases} . $$
We choose $ \delta $ such that the sector defined by $ - 2 \theta_0 - 2 \delta \leq \mbox{arg}(Z) \leq - 2 \theta_0 + 2 \delta $
doesn't meet the sector $ - 2 \theta \leq \mbox{arg}(Z) \leq \epsilon $. We then set
$$ F (Z) = |Z| \exp \Big(  - 2 i \Theta(\mbox{arg}(Z)) \psi (|Z|) + i (1-\psi(|Z|)) \mbox{arg}(Z) \Big) . $$
It is clearly smooth where $ \mbox{arg}(Z) $ is defined hence in the sector $ e^{i(- 2
\theta_0 - 4 \epsilon_1^{\prime} , \epsilon^{\prime} )} ( 0 , +
\infty ) $.  We have $ F (Z) = Z $ for
for $ |Z| \leq C^{-1} $ and  $ |Z| \geq C $ so $ F $ is smooth near $ 0 $. Since
$\Theta(\mbox{arg}(Z)) = \mbox{arg}(Z) $  if $ \mbox{arg}(Z) $ is close to $ -2\theta_0 $, we have (\ref{donnesupport}).
Furthermore, for $
Z $ in the right hand side of (\ref{ensemblededroite}),
we have $ F (Z) - z \ne 0 $ otherwise we should have $ |z| = |Z|
\in [r_1,r_2] $ and then $ z = F (Z) = |z| \exp ( -2 i \Theta(\mbox{arg}(Z))) $
which is impossible by the choice of $ \delta $. This is sufficient to prove (\ref{donneinversible}) since
  $ | F (Z) | \rightarrow \infty $ as $ |Z| \rightarrow \infty
$. \finpreuve

\bigskip

\noindent {\it Proof of Proposition \ref{sjos2}.} We choose first
$ \epsilon_1 $ small enough and $ R_1 $ large enough to ensure
that (\ref{controleellipticite}) and (\ref{controleargument})
hold. We also assume that $   \epsilon_1  $ satisfies the
condition of Lemma \ref{lemmeangulaire}. The full Weyl symbol of 
$ H_0
(\theta_0) + \varepsilon V (\theta_0) $ is of the form
$$ p_{\varepsilon,\theta_0}(x,\xi,h) + h b_{\varepsilon,\theta_0}(x,\xi,h) $$
with $ b_{\varepsilon,\theta_0} $ polynomial of degree $1$ in $ \xi $, and with
\begin{eqnarray*}
 p_{\varepsilon,\theta_0}(x,\xi,h) & = & p_{\varepsilon}^{\rm cl}\left( \kappa_{\theta_0}(x),
 ^t \! \kappa^{\prime}_{\theta_0}(x) \xi,h \right) + a_{\varepsilon}\left( \kappa_{\theta_0}(x),
 ^t \! \kappa^{\prime}_{\theta_0}(x) 
  \xi,h \right)  , \\
 & = : & p_{\varepsilon,\theta_0}^{\rm cl} (x, \xi) + a_{\varepsilon}\left( \kappa_{\theta_0}(x),
 ^t \! \kappa^{\prime}_{\theta_0}(x)  \xi,h \right),
\end{eqnarray*}
 where $ p_{\varepsilon}^{\rm cl} $ is the classical principal symbol and
 $ a_{\varepsilon}(.,.,h) $ a polynomial of degree $1$ in $ \xi $ with coefficients in $ {\mathcal
C}_{\rho}(\theta_0, R_0, \epsilon_0) $, bounded with respect to $h \in (0,h_0 ] $ and $ \varepsilon \in [0,1] $. 
All these symbols are affine (hence smooth) with respect to $ \varepsilon $.
 We then claim that, by possibly increasing $ R_1 $, we may also assume that
 \begin{eqnarray}
  p_{\varepsilon,\theta_0}(x,\xi,h) \in D_F  , \label{conditioncomposition}
 \end{eqnarray}
for all $ h \ll 1 $, $ (x,\xi) \in \Ra^{2d}$ and $ \varepsilon \in
[0,1] $. Note first that, with no loss of generality in Lemma \ref{lemmeangulaire}, we may assume that $ D_F $ is constructed for $ \pi / 2 < \theta_0 < \pi $
so  that $ D_F $ is also a neighborhood of $ \Ra $. Then, for $ |x| \leq R_1 $, $ p_{\varepsilon,\theta_0}(x,\xi,h) $ is real hence belongs to
$ D_F $. On the other hand, there exists $ C_V $ such that 
$$  | a_{\varepsilon}\left( \kappa_{\theta_0}(x),
 ^t \! \kappa^{\prime}_{\theta_0}(x)  \xi,h \right) | \leq C_V R_1^{- \rho} \scal{\xi} , $$
for all  $ R_1 \gg 1 $, $ |x| \geq R_1 $, $ \xi \in \Ra^d $, $ h \in (0,h_0] $ and $ \varepsilon \in [0,1] $. Thus, using (\ref{controleargument}) with $ p_{\iota,\theta_0}^{\rm cl} = p_{\varepsilon,\theta_0}^{\rm cl}  $, we see that for any neighborhood of $ e^{i[-2 \theta_0 - 4 \epsilon_1,\epsilon]}[0,+\infty)  $, we can choose $ R_1 $ large enough such that $ p_{\varepsilon,\theta_0}(x,\xi,h) $ belongs to this neighborhood for $ |x| \geq R_1 $. This implies (\ref{conditioncomposition}) which then shows that
% Denote by $
%p_{\varepsilon,\theta_0}(x,\xi,h) $ the full Weyl symbol of $ H_0
%(\theta_0) + \varepsilon V (\theta_0) $. Then, by
%(\ref{controleargument}), there exists $ C > 0 $ such that,
%\begin{eqnarray}
%  p_{\varepsilon,\theta_0}(x,\xi,h) \in e^{i[- 2 \theta_0 - 4 \epsilon_1 , \epsilon ]} ( 0 , + \infty
%) \cup \{ z \ ; \ |z| \leq C  \} , \nonumber
%\end{eqnarray}
%for all $ h \ll 1 $, $ (x,\xi) \in \Ra^{2d}$ and $ \varepsilon \in
%[0,1] $. 
 $ F \circ p_{\varepsilon,\theta_0} $ is smooth on $ \Ra^{2d} $.  Actually, we have
\begin{eqnarray}
 \psi_{\varepsilon,\theta_0}:= F (  p_{\varepsilon,\theta_0} ) - p_{\varepsilon,\theta_0} \in C_0^{\infty}(\Ra^{2d}) , \label{compactcompact} 
\end{eqnarray}
and, more precisely, $  \psi_{\varepsilon,\theta_0} $ is bounded in $ C_0^{\infty} $ as $ \varepsilon $ and
$h$ vary. Indeed, by (\ref{controleellipticite}), $ |
p_{\varepsilon,\theta_0}(x,\xi,h)| \rightarrow \infty $ as $ |\xi|
\rightarrow \infty $ and, on the other hand, for $ \xi $ in a compact
set, $ p_{\varepsilon,\theta_0}(x,\xi,h) \rightarrow e^{-2 i
\theta_0}|\xi|^2$ as $ |x | \rightarrow \infty $. Using
(\ref{donnesupport}), this gives (\ref{compactcompact}).

 To construct $ K_{\varepsilon} (\theta_0)  $, we recall
the following point. For all $ \Psi \in C_0^{\infty} (\Ra^{2d}) $, we
may write
$$ O \! p_h^w (\Psi) = K (h) + R (h) ,  $$
with $ K (h) $ of finite rank, $ \mbox{rank}(K(h)) \lesssim h^{-d}
 $, and for all $ N \geq 0 $,
$$  || R (h)   ||_{ H^{-N,-N}_{\rm sc}\rightarrow
H^{N,N}_{\rm sc}}  \leq C h^N , \qquad h \ll 1 . $$ In addition, for some
fixed $ \chi \in C_0^{\infty} (\Ra^d) $,
$$ K (h) = \chi K (h) \chi . $$
Let us now
choose $ \Psi \in C_0^{\infty}(\Ra^{2d}) $ such that $ \Psi \equiv 1 $ near
a compact set (independent of $h$ and $ \varepsilon $) containing
the support of $ \psi_{\varepsilon,\theta_0} $.  We then have
$$ O \! p_h^w (\psi_{\varepsilon,\theta_0})  = K (h) O \! p_h^w (\psi_{\varepsilon,\theta_0}) K (h) +
R_{\varepsilon,\theta_0}(h) $$ with, for all $ N \geq 0 $,
$$  || R_{\varepsilon,\theta_0} (h)   ||_{H^{-N,-N}_{\rm sc}\rightarrow H^{N,N}_{\rm sc}}  \leq C h^N ,
\qquad h \ll 1 , \ \varepsilon \in [0,1] , $$ using that $  O \!
p_h^w (\psi_{\varepsilon,\theta_0}) = O \! p_h^w (\Psi)O \! p_h^w
(\psi_{\varepsilon,\theta_0}) O \! p_h^w (\Psi) + {\mathcal
O}(h^{\infty}) $ by pseudodifferential calculus. We then set
$$  K_{\varepsilon} (\theta_0) := K (h) O \! p_h^w (\psi_{\varepsilon,\theta_0}) K (h) . $$
It satisfies (\ref{rank}), (\ref{domaines}) and has a Schwartz
kernel supported in a fixed compact set. To get
(\ref{pasdeparametrix}),  we simply observe that
$$ H_0 (\theta_0) + \varepsilon V (\theta_0)
+ K_{\varepsilon} (\theta_0) - z = O \! p_h^w ( F (
p_{\varepsilon,\theta_0}) - z ) + h T_{\varepsilon}(\theta_0) ,
$$
with $ || T_{\varepsilon}(\theta_0) ||_{H_{\rm sc}^{2,0}
\rightarrow L^2} \lesssim 1 $ as $ h \ll 1 $ and $ \varepsilon \in
[0,1] $. By (\ref{donneinversible}),  $ O \! p_h^w ( F (
p_{\varepsilon,\theta_0}) - z ) $ is invertible for $h$ small
enough (uniformly with respect to $ \varepsilon $ and $ z \in
\Omega $) and so is $ O \! p_h^w ( F ( p_{\varepsilon,\theta_0}) -
z ) + h T_{\varepsilon}(\theta_0) $ by an elementary perturbation
argument. \finpreuve

\medskip

Using the notation of Sj\"ostrand-Zworski \cite{SjostrandZworski}, we now set
\begin{eqnarray}
 \widehat{ H_{\varepsilon}(\theta_0)} =  H_0(\theta_0)  + \varepsilon V (\theta_0)  +
K_{\varepsilon} (\theta_0 ),
  \label{ajoutcompact}
\end{eqnarray}
and
\begin{equation}
 \widetilde{K}_{\varepsilon} (\theta_0 ,z ) = - K_{\varepsilon} (\theta_0 ) ( \widehat{H_{\varepsilon} (\theta_0)}
 -z)^{-1} ,\label{defKtilde}
 \end{equation}
or, equivalently,
\begin{eqnarray}
1+ \widetilde{K}_{\varepsilon} (\theta_0 ,z) = ( H_0 (\theta_0) +
\varepsilon V (\theta_0) - z ) ( \widehat{H_{\varepsilon}
(\theta_0)} - z )^{-1} \label{brancherestelog}
\end{eqnarray}
for all  $z \in  \Omega \setminus {\rm Res}(H_0+\varepsilon
V,\Omega)$. We then have (see \cite{SjBook})
\begin{eqnarray}
\mbox{tr} \left(( H_0 (\theta_0) + \varepsilon
V(\theta_0)-z)^{-1} - ( \widehat{
H_{\varepsilon} (\theta_0) } -z)^{-1} \right) & = &- \mbox{tr}
\left( (1 + \widetilde{K}_{\varepsilon} (\theta_0,z))^{-1}
\partial_z \widetilde{K}_{\varepsilon} (\theta_0, z ) \right)
% \label{tilK0}
\nonumber \\
& = & - \partial_z \log \ \mbox{det}_1 \left( 1 +
\widetilde{K}_{\varepsilon} (\theta_0,z) \right).
\label{deriveelogarithmique}
\end{eqnarray}
Remark that the zeroes of  $ \mbox{det}_1 ( 1 +
\widetilde{K}_{\varepsilon} (\theta_0 , z) ) $ are contained in
the set of resonances since, if $z$ is not a resonance,
(\ref{brancherestelog}) is invertible. Actually,
 the zeroes of $ \mbox{det}_1 ( 1 +
\widetilde{K}_{\varepsilon} (\theta_0 , z) ) $ in $\Omega$ are
exactly the resonances of $H_0+\varepsilon V$ in $\Omega$ with the
same multiplicities (see Definition \ref{defimultiplicite}). More
precisely we recall the following result (see \cite{SjBook}).
\begin{prop}\label{multi}
 If $w \in {\rm Res}(H_0+V, \Omega)$, there exists a holomorphic function $G_w (z)$,
 for $z$ close to $w$, such that $G_w(w) \neq 0$ and
\begin{equation} \label{krein0}
\emph{det}_1\left( 1 + \widetilde{K}_1 (\theta_0 , z )\right) =
(z-w)^{m(w)}G_w(z),
\end{equation}
 where $m(w)$ is the multiplicity of the resonance.
\end{prop}

\noindent {\it Proof.} Let $l(w)$ be the multiplicity of $w$ as
zero of $\hbox{det}_1\left( 1 +
\widetilde{K}_{1}(\theta_0,z)\right)$ given by
\begin{equation}\label{lz0}
l(w)= \frac{1}{2i\pi} \int_{\gamma} \partial_z \log
\hbox{det}_1\left( 1 + \widetilde{K}_1 (\theta_0,z)\right) dz,
\end{equation}
with $\gamma$ a small positively oriented circle centered at $w$.
%Further, we have
%$$ \partial_z \log \hbox{det}_1\Big(1 + K(z)\Big)= \mbox{tr}\Big( (1+K(z))^{-1}\partial_z K(z) \Big),$$
%for any $z \rightarrow K(z)$, operator-valued holomorphic function from $\Omega$ to the trace class.
According to (\ref{deriveelogarithmique}), we have
\begin{eqnarray}
l(w)& =& \frac{i}{2\pi} \int_{\gamma} \mbox{tr} \left(( H_0
(\theta_0) +   V(\theta_0)-z)^{-1}K_1 (\theta_0 ) ( \widehat{ H_1
(\theta_0) }  -z)^{-1}
\right)dz \nonumber \\
&  = &\frac{i}{2\pi} \mbox{tr} \left(  \int_{\gamma} ( H_0
(\theta_0) + V(\theta_0)-z)^{-1} - ( \widehat{ H_1 (\theta_0)  }
-z)^{-1} dz \right)  \nonumber.
\end{eqnarray}
By construction of $ \widehat{ H_1 (\theta_0)  }$, the resolvent $ (
\widehat{ H_1 (\theta_0)  } -z)^{-1}$ is holomorphic near $w$ and
its integral on $\gamma$ vanishes. It follows that
$$l(w)=  \mbox{tr} \left( \frac{i}{2\pi} \int_{\gamma} \
( H_0 (\theta_0) +   V(\theta_0)-z)^{-1} dz \right) = \mbox{tr}(
\Pi_{\theta_0,w}), $$ where $ \Pi_{\theta_0,w}$, defined by
(\ref{defprojecteur}),  is a projector which (by definition of the
multiplicity $m(w)$) satisfies
$$\mbox{tr} ( \Pi_{\theta_0,w}) = {\rm rank} ( \Pi_{\theta_0,w})= m(w).$$
This conclude the proof of Proposition \ref{multi}. \finpreuve

\bigskip

Therefore, the multiplicities of the resonances as zeroes of $
\hbox{det}_1\left( 1 + \widetilde{K}_1 (\theta_0 ,z)\right) $ or
as given by Definition \ref{defimultiplicite} coincide and we have
the factorization
\begin{eqnarray}
 \hbox{det}_1\left( 1 + \widetilde{K}_1 (\theta_0 , z )\right) = \prod_{w \in
    \rm{Res}(H_0+V,\Omega)} (z-w) G_1 (z,h) \label{factorisationSjostrand}
\end{eqnarray}
where, for each $ h \in (0,h_1] $, $ G_1 (.,h) $ is a non
vanishing holomorphic function on $ \Omega $. 

We now recall a beautiful
result due to Sj\"ostrand which is a crucial consequence of Proposition  \ref{sjos2}.
\begin{prop}[\cite{Sjosequi}] \label{Harnack} There exists $ \varphi^{G_1} (.,h) $ holomorphic on $ \Omega $ such that
$$ G_1 (z,h) = \exp \left( \varphi^{G_1} (z,h) \right), \qquad  \ h \ll 1, \ z \in \Omega, $$
and, for all $ W \Subset \Omega $
$$ | \partial_z \varphi^{G_1} (z,h) | \leq C_W h^{-d} , \qquad h \ll 1 , \ z \in W .$$
\end{prop}
An immediate consequence of (\ref{factorisationSjostrand}) and
Proposition \ref{Harnack} is that, for all $ W \Subset \Omega $,
\begin{eqnarray}
\left| \partial_z^k \log  \hbox{det}_1 \left( 1 + \widetilde{K}_1
(\theta_0,z)\right) - \sum_{w \in
    \rm{Res}(H_0+V,\Omega)} \frac{(k-1)!}{(w-z)^k} \right| \leq C_W h^{-d}, \label{resonanceszeta1}
\end{eqnarray}
for $ h \ll 1 $ and $ z \in W $.
 The same result applied with $ V \equiv 0 $, using that $ H_0
$ has no resonances, shows that
\begin{eqnarray}
\left| \partial_z^k \log  \hbox{det}_1 \left( 1 +
\widetilde{K}_{0}(\theta_0,z)\right) \right| \leq C_W h^{-d} ,
\label{resonanceszeta0}
\end{eqnarray}
for $ h \ll 1 $ and $ z \in W $.

Another useful consequence of the absence of resonance for $ H_0$
is the following. Since $ H_0 $ has no resonances, $ H_0 (\theta_0)
- z  $ is invertible for all  $ h \ll 1 $ and  all $z  $ in a
neighborhood of $ \overline{\Omega}  $.
 Therefore, for all $ h \ll 1 $,
there exists $ \varepsilon_h   $,  such that $   H_0 (\theta_0) +
\varepsilon V(\theta_0)-z  $ is invertible for $ | \varepsilon | <
\varepsilon_h   $ and $ z \in \Omega $. Thus,  by
(\ref{brancherestelog}), the function
\begin{equation}\label{defG}
G_{\varepsilon}(z,h):= \mbox{det}_1 \left( 1 +
  \widetilde{K}_{\varepsilon} (\theta_0,z) \right),
\qquad z \in \Omega , \ \varepsilon \in (- \varepsilon_h ,
\varepsilon_h),
\end{equation}
is holomorphic and doesn't vanish. This allows to choose a branch
of its logarithm which we denote by $ \mbox{Log}_h \
G_{\varepsilon}(z,h)  $,  to stress on the $h$ dependence of such
a choice.

\begin{prop} \label{loginitial} The branch $ \emph{Log}_h \  G_{\varepsilon}(z,h) $ can
be chosen such that,  given  a fixed $ z_0 \in \Omega^+_{\delta} $, we
have, for all $ j \geq 0 $, $ l \geq 1 $,
$$ \left| \frac{\partial^l}{\partial z^l}  \frac{\partial^j}{\partial \varepsilon^j} \emph{Log}_h \  G_{\varepsilon}
(z_0,h )_{| \varepsilon = 0} \right| \lesssim h^{-d} .
$$
\end{prop}

\noindent {\it Proof.} According to (\ref{resonanceszeta0}),  $ G_0 (z,h) = \exp
(\varphi^{G_0}(z,h)) $ with $ | \partial_z \varphi^{G_0} (z,h) | \lesssim
h^{-d} $ . On the other hand, for all $ h \ll 1 $, we can find $
\varepsilon (z_0,h)
> 0 $ such that
$$ \left|  \frac{G_{\varepsilon} (z_0,h)}{G_0(z_0,h)} - 1  \right| \leq 1/2, \qquad |\varepsilon |
\leq \varepsilon (z_0,h)  $$ thus we can set
\begin{eqnarray}
\mbox{Log}_h \  G_{\varepsilon} (z_0,h) = \varphi^{G_0} (z_0,h) + \log
\left( \frac{G_{\varepsilon}(z_0,h)}{G_0(z_0,h)} \right)
\label{conditionenz0}
\end{eqnarray}
 where $ \log $ is
the principal determination of the logarithm on $ \Ca \setminus (-
\infty , 0 ] $. We can then define  $ \mbox{Log}_h \
G_{\varepsilon}(z,h) $ as the unique primitive of $
\partial_z G_{\varepsilon}(z,h) / G_{\varepsilon}(z,h) $
 coinciding with the right hand side of (\ref{conditionenz0})
 at $ z=z_0 $. The smoothness  with respect to $z$ and $ \varepsilon
 $ (close to $0$) is then clear. The bounds on
 the derivatives at $ z = z_0 $ and $ \varepsilon = 0 $ are
 obtained by applying $ \partial_{\varepsilon}^k \partial_z^{l-1} $
to (\ref{deriveelogarithmique}), using Proposition \ref{sjos2} and
(\ref{borneresolventedelta}).
%\begin{eqnarray}
% \frac{\partial}{ \partial \varepsilon} \mbox{Log}_h \ G_{\varepsilon}(z,h) = \mbox{tr} \left(
%(1+\widetilde{K}_{\varepsilon} (\theta_0,z))^{-1} \frac{\partial}{
%\partial \varepsilon} \widetilde{K}_{\varepsilon} (\theta_0,z)
% \right) . \nonumber
%\end{eqnarray}
%Using (\ref{brancherestelog}), Propositions
%\ref{Omegadeltaresolvente} and \ref{sjos2}, we have
%$$ || (1+\widetilde{K}_0 (\theta_0,z_0))^{-1}||_{L^2 \rightarrow L^2} \lesssim 1,
%\qquad
%  \left\| \ \frac{\partial}{ \partial \varepsilon}
%  \widetilde{K}_{\varepsilon}
%  (\theta_0,z_0)_{| \varepsilon = 0} \right\|_{\rm tr} \lesssim h^{-d} , $$ which yield the
%bound $ |\partial_{\varepsilon} \mbox{Log}_h \ G_{\varepsilon}
%(z_0,h)_{| \varepsilon = 0} | \lesssim h^{-d} $. We obtain the
%same bound for higher order derivatives by induction.
\finpreuve

\bigskip

Regarding the behavior of $ \partial_{\varepsilon}^j \mbox{Log}_h
\ G_{\varepsilon} (z,h)_{|\epsilon = 0} $ for $ z \in \Omega $, we
have the following result.

\begin{prop} \label{pourtheoremeavecborne} For all $ j \geq 0 $, $ l \geq 1 $, there exists $ N_{j,l} \in \Na $ such that, for  all compact
subset $ W \Subset \Omega $,
$$ \left| \partial_{\epsilon}^j \partial_z^l \emph{Log}_h
\ G_{\varepsilon} (z,h)_{|\epsilon = 0} \right|\leq C_{W} h^{-d} \sup_{Z \in \Omega}
\left( 1 + || (H_0(\theta_0)-Z)^{-1}||_{L^2 \rightarrow H^{2,0}_{\rm
sc }} \right)^{N_{j,l}} , \qquad h \ll 1, \ z \in W . $$
\end{prop}

\noindent {\it Proof.} %Consider first $ j  = 0 $ and $ l = 1 $. 
By
 writing $ \mbox{Log}_h \
G_{\varepsilon} (z,h)_{|\epsilon = 0} $ as the sum of $
\mbox{Log}_h \ G_{\varepsilon} (z_0,h)_{|\epsilon = 0} $ and the
integral of its derivative over a path joining $z_0$ to $z$, the
result follows from (\ref{deriveelogarithmique}), (\ref{defG}),  Proposition \ref{sjos2}.
and Proposition \ref{loginitial}.
\finpreuve

\section{Proofs of Theorems  \ref{theoremeprincipalp}, \ref{theoremeprincipal2} and
\ref{theoremeavecborne}}
\setcounter{equation}{0}
% In this section, we prove the existence of a meromorphic
%continuation (with respect to $z$) of the Zeta function  $ \zeta_p (k,z,h)  $ for $ k $ large enough. Using Proposition
%\ref{reductionZetaintro}, this will immediately imply Theorems
%\ref{theoremeprincipalp}  and \ref{theoremeprincipal2}.
%We shall
%also prove in passing a slightly more explicit factorization of
%the regularized determinant involving the operators $
%\widetilde{K}_{\varepsilon}(\theta_0,z)  $ introduced in
%subsection \ref{pourlafactorisation}.

\subsection{The general case} \label{preuvedutheoreme}
 Using the notation  (\ref{Taylorkp}), we have, for $ k  > d/2 $,
\begin{eqnarray}
 \zeta_{p} (k,z,h) = T_p^k (H_0,V,z), \qquad   h \ll 1 , \ \  z \in \Omega^+  , \label{pointdedepart}
\end{eqnarray}
and, by Proposition \ref{deformationSjp}, we also have, if $ k  > d/2 +  1 $ which we now assume,
\begin{eqnarray}
 T_p^k (H_0,V,z) =  T_p^k (H_0(\theta_0),V(\theta_0),z)  , \qquad
 h \ll 1 , \ \  z \in \Omega^+  . \label{dilatationSj}
\end{eqnarray}
To analyze the right hand side of (\ref{dilatationSj}), we consider first 
$$ \widehat{T}_p^k (\theta_0,z,h) := \emph{tr} \left( ( \widehat{ H_{1} (\theta_0)  }
-z)^{-k} - \sum_{j=0}^{p-1} \frac{1}{j!} \frac{d^j}{d
\varepsilon^j} {( \widehat{ H_{\varepsilon} (\theta_0)  } -z)^{-k}}_{ | \varepsilon = 0}
\right) ,
$$
where $ \widehat{ H_{\varepsilon} (\theta_0)  } $ is defined by (\ref{ajoutcompact}).
 % By (\ref{ajoutcompact}), we have the resolvent identity
%\begin{eqnarray}
%( H_0 (\theta_0) + \varepsilon  V(\theta_0)-z)^{-1}= ( \widehat{ H_{\varepsilon} (\theta_0)  }
%-z)^{-1} + ( H_0 (\theta_0) + \varepsilon V(\theta_0)-z)^{-1}
%K_{\varepsilon} (\theta_0 ) ( \widehat{ H_{\varepsilon} (\theta_0)
%} -z)^{-1} , \label{hatresolv}
%\end{eqnarray}
 %at least for $ h \in (0,h_0]  $, $ z \in \Omega^+  $ and $ \varepsilon  = 1  $ or close
%to $ 0 $.
\begin{lemm} \label{lemmenoncompactetrace} For all $ h \ll 1 $, the function $ \widehat{T}_p^k (\theta_0,z,h) $
is well defined, has an holomorphic continuation from $ \Omega^+ $ to  $   \Omega$
and, for all $ W \Subset \Omega $,
$$ | \widehat{T}_p^k (\theta_0,z,h) | \leq C_{W} h^{-d} , \qquad h \ll 1 , \ z \in W . $$
\end{lemm}

\noindent {\it Proof.} Write first that
\begin{eqnarray}
 \frac{d}{d \varepsilon} (\widehat{ H_{\varepsilon} (\theta_0)  }
-z)^{-1} = -
 (\widehat{ H_{\varepsilon} (\theta_0)  } -z)^{-1} \left( V(\theta_0) + \partial_{\epsilon} K_{\varepsilon}
 (\theta_0) \right) ( \widehat{ H_{\varepsilon} (\theta_0)  } -z)^{-1}. \label{deriveeenepsilon}
\end{eqnarray}
Then, an elementary induction shows that the operator
$$\frac{d^j}{d
\varepsilon^j} ( \widehat{ H_{\varepsilon} (\theta_0)  } -z)^{-1}- j! \left( -(
\widehat{ H_{\varepsilon} (\theta_0) } -z)^{-1} V(\theta_0) \right)^j(
\widehat{ H_{\varepsilon} (\theta_0) } -z)^{-1}$$
is a linear combination of holomorphic finite rank operators with trace norm of order $ h^{-d}$, for all $ j $.
This  formula for $ j =p$ combined with Taylor's formula and Proposition \ref{sjos2} shows that the operator
$$  ( \widehat{ H_{1} (\theta_0)  }
-z)^{-1} - \sum_{j=0}^{p-1} \frac{1}{j!} \frac{d^j}{d
\varepsilon^j} ( \widehat{ H_{\varepsilon} (\theta_0)  } -z)^{-1}_ { | \varepsilon = 0}+ p
\!
 \int_0^1 (\varepsilon-1)^{p-1} \left( (
\widehat{ H_{\varepsilon} (\theta_0) } -z)^{-1} V(\theta_0) \right)^p(
\widehat{ H_{\varepsilon} (\theta_0 } -z)^{-1} d \varepsilon
$$
is a linear combination of holomorphic trace class operators with
norm $ {\mathcal O}(h^{-d}) $, locally uniformly on compact subsets of $ \Omega $. Using (\ref{pourderiverresolvente}), Proposition \ref{pourlimite} and (\ref{pasdeparametrix}), the
$ k $-th derivative of the operator in the integral above is trace class,
holomorphic on $ \Omega $ and with trace norm $ {\mathcal O}(h^{-d}) $, locally uniformly with respect to $z$. The result follows.
\finpreuve

\bigskip

Using (\ref{pourderiverresolvente}) and  (\ref{dilatationSj}), we
obtain
\begin{eqnarray}
 T_p^k (H_0,V,z) =  \widehat{T}_p^k \left( \theta_0,z,h \right) + \frac{1}{(k-1)!} \partial_z^{k-1}
 A(z,h) , \qquad h \ll 1 , \ z \in \Omega^+ , \label{pourfacteurderive}
\end{eqnarray}
 where
\begin{eqnarray}
A (z,h) & = & \mbox{tr} \left( ( H_0 (\theta_0) +
V(\theta_0)-z)^{-1} - ( \widehat{ H_1
(\theta_0) } -z)^{-1}  \right)  \nonumber \\
& & \ \ \ \ - \sum_{j=0}^{p-1} \frac{1}{j!} \frac{d^j}{d
\varepsilon^j} \mbox{tr} \left(  ( H_0 (\theta_0) + \varepsilon
V(\theta_0)-z)^{-1} -  ( \widehat{
H_{\varepsilon} (\theta_0) } -z)^{-1} \right)_{ | \varepsilon = 0} ,
\nonumber
\end{eqnarray}
that is
 \begin{eqnarray}
 - A (z,h) =  \partial_z \log \ \mbox{det}_1 \left( 1 +
  \widetilde{K}_1 (\theta_0,z) \right) - \partial_z
  \sum_{j=0}^{p-1} \frac{1}{j!}\frac{d^j}{d \varepsilon^j}
\ \mbox{Log}_h  {G_\varepsilon (z,h)}_{| \varepsilon = 0} .
\label{autiliserbis}
\end{eqnarray}
 by (\ref{deriveelogarithmique}), (\ref{defG}) and the notation of Propositions \ref{loginitial} and \ref{pourtheoremeavecborne}.
\bigskip

\noindent {\bf Proof of Theorems \ref{theoremeprincipalp} and
\ref{theoremeavecborne}.} By (\ref{pointdedepart}),
(\ref{dilatationSj}), (\ref{pourfacteurderive}),
(\ref{autiliserbis}) and (\ref{factorisationSjostrand})
 we have an expression of the form
(\ref{aintegrerkfois}) with
\begin{equation}
 \phi_p (z,h) = \widehat{T}_p^k \left( \theta_0,z,h \right) - \frac{1}{(k-1)!}  \partial_z^{k} \Big(  \varphi^{G_1}(z,h) -
  \sum_{j=0}^{p-1} \frac{1}{j!}\frac{d^j}{d \varepsilon^j}
\ \mbox{Log}_h  {G_\varepsilon (z,h)}_{| \varepsilon = 0} \Big)
\label{expphip}
\end{equation}
which is holomorphic on $
\Omega $. This proves Theorem \ref{theoremeprincipalp} using
Proposition \ref{reductionZetaintro} with $ {\mathcal
H}(\Omega,h_1) $ the set of families of holomorphic functions on $
\Omega $.

To prove Theorem \ref{theoremeavecborne}, we simply additionally
note that, by Proposition \ref{Harnack} and Proposition
\ref{pourtheoremeavecborne}, we can find $ N > 0 $ such that, for all $ W \Subset
\Omega $,
\begin{eqnarray}
 |\phi_p (z,h)| \leq C_W h^{-d} \sup_{Z \in \Omega} \left( 1 + ||
(H_0(\theta_0)-Z)^{-1}||_{L^2 \rightarrow H^{2,0}_{\rm sc }}
\right)^{N}  , \qquad h \ll 1, \ z \in W .
%\label{borneespaceholomorphe}
\end{eqnarray}
  Then, Proposition \ref{reductionZetaintro} gives the
result using the space $ {\mathcal H} (\Omega,h_1) $ of families
of holomorphic functions locally bounded by (a constant times) $ h^{-d} \sup_{Z \in \Omega} \left( 1 + ||
(H_0(\theta_0)-Z)^{-1}||_{L^2 \rightarrow H^{2,0}_{\rm sc }}
\right)^{N} $.
% the
%right hand side of (\ref{borneespaceholomorphe}). 
Note that it
satisfies (\ref{stabiliteconstantes}) and (\ref{stabiliteparprimitive}). \finpreuve

\subsection{Proof of Theorem \ref{theoremeprincipal2}} \label{preuvedutheoreme2}
In this subsection, $ {\mathcal H} (\Omega,h_1) $ denotes the space of
families of holomorphic functions $ (\phi (.,h))_{h \in (0,h_1]} $
such that, for all $ W \Subset \Omega $, $ |\phi (z,h)| \leq C_{W}
h^{-d} $, for $z \in W $ and $ h \in (0,h_1] $.

\medskip

 For $ p = 1 $, the result can be considered as essentially a
consequence of \cite{BruneauPetkov}. For completeness, we give the
proof. In that case, $\phi_1$ (given by (\ref{expphip}) with $p=1$) belongs to $ {\mathcal H} (\Omega,h_1) $ according to Lemma
\ref{lemmenoncompactetrace}, Proposition \ref{Harnack} and  (\ref{resonanceszeta0}).   The result follows then from  Proposition
\ref{reductionZetaintro}.
%using Proposition \ref{Harnack} (also for $ V \equiv 0 $), (\ref{autiliserbis}) is reduced to
%\begin{eqnarray}
% \partial_z \log \mbox{det}_1 \left( 1 +
%  \widetilde{K}_{1} (\theta_0,z) \right) - \partial_z \log \mbox{det}_1 \left( 1 +
%  \widetilde{K}_0 (\theta_0,z) \right) = \sum_{w \in {\rm Res}(H_0+V,\Omega)}
%\frac{1}{z-w} + \phi(z,h) \label{contributionresonances}
%\end{eqnarray}
%with $ \phi \in {\mathcal H} (\Omega,h_1) $. Then, using Lemma
%\ref{lemmenoncompactetrace} and Proposition
%\ref{reductionZetaintro}, the result follows.

\medskip

 In the case $ p = 2 $,  (\ref{expphip}) gives
 %the term $  A (z,h) $  is the the sum of $ - \partial_z^{k-1}
%(\ref{contributionresonances})/ (k-1)! $ and
$$\phi_2(z,h)- \widehat{T}_2^k \left( \theta_0,z,h \right)+ \frac{1}{(k-1)!}  \partial_z^{k}  \varphi^{G_1}(z,h)=$$
\begin{eqnarray}
\frac{1}{(k-1)!} \partial_z^{k}
\ \mbox{Log}_h  {G_0 (z,h)}_{| \varepsilon = 0} +
\mbox{tr} \left( \frac{d}{d \varepsilon}
(\widehat{H_{\varepsilon}(\theta_0)}-z)^{-k}_{ | \varepsilon = 0} -
\frac{d}{d \varepsilon}
  (H_0(\theta_0)+ \varepsilon
V(\theta_0)-z)^{-k}_{| \varepsilon = 0 } \right) . \label{relaphi2}
\end{eqnarray}
By Lemma \ref{lemmenoncompactetrace}, Proposition \ref{Harnack} and (\ref{resonanceszeta0}),  it remains to study
 the second term of (\ref{relaphi2}). We first remark that this term can be
 written as the sum of
 \begin{eqnarray}
   \mbox{tr} \left( \frac{d}{d \varepsilon}
(\widehat{H_{0}(\theta_0)} + \varepsilon V (\theta_0)-z)^{-k}_{ |
\varepsilon = 0} - \frac{d}{d \varepsilon}
  (H_0(\theta_0)+ \varepsilon
V(\theta_0)-z)^{-k}_{| \varepsilon = 0 } \right)
\label{atransformer}
\end{eqnarray}
 and
$$ - \partial_z^{k-1} \mbox{tr} \left( (\widehat{H_{0}(\theta_0)}
-z)^{-1}\partial_{\varepsilon} K_{\varepsilon} (\theta_0)_{|
\varepsilon = 0}
 (\widehat{H_{0}(\theta_0)} -z)^{-1} \right) / (k-1)!, $$
using (\ref{ajoutcompact}) and (\ref{deriveeenepsilon}). This last expression  clearly
belongs to $ {\mathcal H}(\Omega,h_1) $ by Proposition \ref{sjos2} and we are left with the
study of (\ref{atransformer}).

 For that purpose, we use
the approximation $ V_n $ of $ V $ introduced in Lemma
\ref{approxdistordu}. Using (\ref{deriveeTaylor}), Lemma
\ref{lemmeclassiquetrace} and an elementary cyclicity argument, we
can write
\begin{eqnarray}
\mbox{tr} \left( \frac{d}{d \varepsilon} (\widehat{H_0(\theta)}+ \varepsilon
  V_n(\theta)-z)^{-k}_{ | \varepsilon = 0} \right) = - k
\mbox{tr} \left(  V_n(\theta)  (\widehat{H_0(\theta)}-z)^{-k-1} \right) . \label{partiesansprobleme}
\end{eqnarray}
Writing  $ \frac{d}{d \varepsilon}
  (H_0(\theta)+ \varepsilon
V_n(\theta)-z)^{-k} $ as the derivative of $ (H_0(\theta)+
\varepsilon V_n(\theta)-z)^{-k} -  (H_0(\theta) -z)^{-k} $ with
respect to $ \varepsilon  $ and using Proposition
\ref{deformationSjtrace}, we obtain similarly
\begin{eqnarray}
\mbox{tr} \left( \frac{d}{d \varepsilon} (H_0(\theta)+ \varepsilon
  V_n(\theta)-z)^{-k}_{ | \varepsilon = 0} \right) = - k
\mbox{tr} \left(  V_n  (H_0-z)^{-k-1} \right) . \label{limiteacalculer}
\end{eqnarray}
Substracting  $ - k \mbox{tr} \left(  V_n(\theta) (e^{-2i\theta}
  H_0-z)^{-k-1}  \right) $ to (\ref{partiesansprobleme}) and
(\ref{limiteacalculer}) and then letting $ n \rightarrow \infty $
using Proposition \ref{lemmeprop}, (\ref{atransformer}) can thus
be written as the sum of
\begin{eqnarray}
- k  \mbox{tr} \left(  V(\theta) \big( (\widehat{H_0(\theta)}
-z)^{-k-1} - (e^{-2i\theta}
  H_0-z)^{-k-1} \big)  \right) \label{wellbehaved}
\end{eqnarray}
and
\begin{eqnarray}
 \lim_{n \rightarrow \infty}  k \mbox{tr} \left(  V_n  (
  H_0-z)^{-k-1}  -    V_n(\theta)  (e^{-2i\theta}
  H_0-z)^{-k-1}  \right) . \label{explicite}
\end{eqnarray}

\begin{prop}\label{prop52}
 (\ref{wellbehaved}) belongs to $ {\mathcal H}(\Omega,h_1) $.
 \end{prop}

\noindent {\it Proof.} By the resolvent identity, (\ref{grandt})
and Proposition \ref{sjos2}, we have
$$ (\widehat{H_0(\theta)}
-z)^{-1} - (e^{-2i\theta}
  H_0-z)^{-1} = (\widehat{H_0(\theta)}
-z)^{-1} B (h) (e^{-2i\theta}
  H_0-z)^{-1} , $$
with $ B (h) = O \! p_h^w (b(h)) $ for some  family $ (b (h))_{h \ll 1} $
bounded in $ S^{2,-N} $ for all $ N $. Using
(\ref{pourderiverresolvente}), the operator $ V(\theta) \big(
(\widehat{H_0(\theta)} -z)^{-k-1} - (e^{-2i\theta}
  H_0-z)^{-k-1} \big) $ is therefore a linear combination of
  operators of the form
  $$  V(\theta)  (\widehat{H_0(\theta)}
-z)^{-k_1-1} B (h) (e^{-2i\theta}
  H_0-z)^{-k_2-1} , \qquad k_1 + k_2 = k  . $$
By (\ref{normrestriction}), (\ref{pasdeparametrix}) and Lemma
\ref{lemmeclassiquetrace}, each such operator has a trace norm of
order $ h^{-d} $, uniformly with respect to $ z \in \Omega $, so
the result follows. \finpreuve

\begin{prop}\label{prop53}
(\ref{explicite}) belongs to $ {\mathcal H}(\Omega,h_1) $.
\end{prop}

\noindent {\it Proof.} The operators $ V_n  (
  H_0-z)^{-k-1}  $ and $   V_n(\theta)  (e^{-2i\theta}
  H_0-z)^{-k-1} $ are both trace class so we compute their traces
  separately.
By writing
$$ V_n (\theta) = \sum_{|\alpha| \leq 2}
  v_{n,\alpha,\theta}(x,h) (h D)^{\alpha} , $$
  we first have
\begin{eqnarray}
\mbox{tr} \left( V_n(\theta)  (e^{-2i\theta}
  H_0-z)^{-k-1} \right) = (2 \pi h)^{-d} \int \! \! \int_{\Ra^{2d}}\sum_{|\alpha| \leq 2}
  v_{n,\alpha,\theta}(x,h) \xi^{\alpha}
  (e^{-2 i \theta} \xi^2 -
  z)^{-k-1} dx d\xi .
\end{eqnarray}
This holds also for $ \theta = 0 $ which gives an expression for $
\mbox{tr} \left( V_n  (
  H_0-z)^{-k-1} \right) $. In the latter case, deforming $ \Ra^d_{\xi} $ into $ e^{-i \theta} \Ra^d_{\xi}
  $, we get
$$ \mbox{tr} \left( V_n  (
  H_0-z)^{-k-1} \right) = (2 \pi h)^{-d} \int \! \! \int_{\Ra^{2d}}\sum_{|\alpha| \leq 2}
  v_{n,\alpha,0}(x,h) ( e^{-i \theta}\xi )^{\alpha}
  (e^{-2 i \theta} \xi^2 -
  z)^{-k-1}  e^{-i d \theta} d\xi dx , $$
and the last integral can be rewritten as
\begin{eqnarray}
(2 \pi h)^{-d} \int \! \! \int_{\Ra^{2d}}\sum_{|\alpha| \leq 2}
  v_{n,\alpha,0}(\kappa_{\theta}(x),h) ( e^{-i \theta}\xi )^{\alpha}
  (e^{-2 i \theta} \xi^2 -
  z)^{-k-1}  e^{-i d \theta} d\xi \ \mbox{det} (\partial_x \kappa_{\theta}(x))
  dx.
\end{eqnarray}
To justify this last deformation, one simply notices that $ \int
v_{n,\alpha,0}(\kappa_{\theta}(x),h) \ \mbox{det} (\partial_x
\kappa_{\theta}(x))
  dx $ depends holomorphically on $ \theta $ and that it is
  constant for $ i \theta $ real and close to zero since $ \kappa_{\theta}
  $ is then a diffeomorphism from $ \Ra^d $ to itself. Now for $ |x| \geq R
  $  large enough, (independent of $n$), we have $ \kappa_{\theta}(x) = e^{i \theta}x
  $ and
  $$  v_{n,\alpha,0}(\kappa_{\theta}(x),h) e^{-i |\alpha| \theta} = v_{n,\alpha,\theta}(x,h),
  \qquad  e^{-i d \theta } \ \mbox{det} (\partial_x \kappa_{\theta}(x))  = 1 . $$
Therefore, if we set
$$ c_{n,\alpha,\theta} (x,h) = v_{n,\alpha,0}(\kappa_{\theta}(x),h) e^{-i |\alpha| \theta}
 e^{-i d \theta } \ \mbox{det} (\partial_x \kappa_{\theta}(x))  - v_{n,\alpha,\theta}(x,h) $$
which is compactly supported, we have
$$ (\ref{explicite}) = \lim_{n \rightarrow \infty} k\sum_{|\alpha| \leq 2} (2 \pi h )^{-d}
 \int_{\Ra^d} \xi^{\alpha}(e^{-2 i \theta}\xi^2 -z)^{-k-1} d \xi \times
 \int_{|x| \leq R} c_{n,\alpha,\theta} (x,h) dx , $$
which is easily seen to belong to $ {\mathcal H}(\Omega,h_1) $. \finpreuve

\bigskip

 The conclusion follows then from (\ref{relaphi2}), Propositions \ref{prop52}, \ref{prop53} and \ref{reductionZetaintro}.
\finpreuve

\section{A counter example for $ p = 3 $} \label{preuvecontreexemple}
\setcounter{equation}{0}
In this section, we prove Theorem \ref{theoremeprincipal3}.
%\ref{theoremeprincipal2} cannot be improved in general, ie that if
%$ p \geq 3 $ there exists $ V \in {\mathcal V}_{\rho} $ for which
%the holomorphic function $ \varphi_p$ doesn't belong to . We give
%a counterexample for $d=1$ and $ p =3 $ with a compactly supported
%potential.
We  consider $ H_0 = - h^2 \frac{d^2}{d x^2} $ on $ L^2
(\Ra) $ and $ V $ a compactly supported bounded potential. In that
case $ V (H_0 - z)^{-1} $ is in the trace class for all $ z \notin [ 0
, + \infty ) $ hence in any Schatten class $ \mathbf{S}_p $.  For
trace class operators $ K \in \mathbf{S}_1 $, the formula
(\ref{detFredholm}) can be written
$$ \mbox{Det}_p (I+K) = \mbox{Det}_1 (I+K) \exp \left( \sum_{j=1}^{p-1} \frac{(-1)^j}{j} \mbox{tr}(K^j)  \right) . $$
We therefore obtain
\begin{eqnarray}
  D_3 (H_0,H_0 + V ; z,h) = D_2 (H_0,H_0 + V ;
z,h ) e^{\frac{1}{2} \phi(z,h)} \label{avecfacteurcorrectif}
\end{eqnarray}
 where
$$ \phi (z,h) = \mbox{tr} \left( V(H_0-z)^{-1}V(H_0-z)^{-1} \right) . $$
For $ z = k^2 $ with $ \mbox{Im}(k) > 0 $, the integral kernel of
$ ( H_0 - z )^{-1} $ is $ i e^{i k |x-x^{\prime}|/h} / (2hk) $ and $
\phi (z,h) $ can be computed explicitly, namely
\begin{eqnarray}
 \phi (k^2,h) & = & \frac{-1}{(2h k)^2} \int_\Ra \int_\Ra V(x)
V(x^{\prime})e ^{2ih^{-1}k | x-x^{\prime} | } dxdx^{\prime} ,
\nonumber \\
&  =  & \frac{-1}{(2h k)^2} \int_\Ra \widetilde{V}(y) e
^{2ih^{-1}k \vert y \vert} dy, \label{a2W}
\end{eqnarray}
with
\begin{eqnarray}
\widetilde{V}(y) = \int_\Ra V(x) V(x-y)dx. \label{defW}
\end{eqnarray}
% Notice that, if $ V $ is even, then $ \widetilde{V} $ is even and $ \widetilde{V} = V \ast V $.
 Setting
$$ \widetilde{V}^{+}_{\rm ev}(y) = {\bf{1}}_{[0,+\infty)} \left( \widetilde{V} (y) + \widetilde{V} (-y) \right)  , $$
we have
\begin{eqnarray}
  \phi (k^2,h) =  - \frac{2 \pi}{(2 k h)^2} ({\mathcal F}_{\rm
inv} \widetilde{V}^+_{\rm ev}) (2 k h^{-1}) ,
\label{definitionphi}
\end{eqnarray}
 where $ {\mathcal
F}_{\rm inv} $ is the usual inverse Fourier transform
$$  {\mathcal F}_{\rm inv} g(\xi) = \frac{1}{2 \pi} \int e^{ ix \xi } g(x) dx . $$
For example, for the characteristic function $V(x) = \chi_a (x) :=
{\bf{1}}_{[-a,a]}(x)$, we have
$$ \widetilde{V} (y) = \begin{cases}  (2a - y)_+ & \mbox{if} \ y \geq 0 , \\
(2a + y)_+ & \mbox{if} \ y < 0 , \end{cases}  $$ where $ (t)_+ =
\max (t,0) $. After elementary computations, we also obtain in
this explicit case
$$ \phi (k^2,h)= \frac{-ia}{2 k^3h}+ \frac{1}{8 k^4}(e ^{4iah^{-1}k}-1). $$
%
%For example for the characteristic function $V(x) =
%{\bf{1}}_{[0,a]}(x)$, we have
%$$ \widetilde{V} (y) = \begin{cases}  (a-y)_+ & \mbox{if} \ y \geq 0 , \\
%(a+y)_+ & \mbox{if} \ y < 0 , \end{cases}  $$ where $ (t)_+ = \max
%(t,0) $. By straightforward computation, we also obtain in this
%explicit case
%$$ \phi (k^2,h)=  \frac{-ia}{4 k^3h}+ \frac{1}{8 k^4}(e ^{2iah^{-1}k}-1) .  $$
For $ k = z^{1/2} $ with $ \mbox{Im} (k) < 0 $, which makes sense
at least close to $ 1 $, this examples shows that
$$ | \partial_z \phi (z,h) | \gtrsim \exp \left(   a | \mbox{Im}(k) | / h  \right) , \qquad h \ll 1 . $$
This proves that the logarithmic derivative of the corrective
factor in (\ref{avecfacteurcorrectif}) can indeed blow up
exponentially, which is a strong form of the estimate (\ref{bornepresqueexponentielle}).

This elementary striking example doesn't however fit in our
framework since $ V $ is not smooth. In particular, it can not be
used directly to prove Theorem \ref{theoremeprincipal3}. For the
latter proof, we need the following lemma.

%Fix $ a = 1 $ and let $ V = \chi_1 $ as above. Choose a family of
%even functions $ V_{\epsilon} \in C_0^{\infty}(\Ra,\Ra) $ such
%that
%$$  0 \leq V_{\epsilon} \leq 1, \qquad \mbox{supp} (V_{\epsilon}) \subset [-1,1] , \qquad  V_{\epsilon}
%\equiv 1 \ \ \mbox{on} \ \ [\epsilon-1,1-\epsilon] . $$ Then $ (
%V_{\epsilon} )_{0 < \epsilon \ll 1} $ is bounded in $ L^{\infty}
%(\Ra) $ and $ V - V_{\epsilon} \rightarrow 0 $ in $ L^1 (\Ra) $ so
%$$ || V_{\epsilon} \ast V_{\epsilon} - V \ast V ||_{L^{\infty}(\Ra)} \rightarrow 0 , \qquad \epsilon \rightarrow 0 . $$
%Since $ \widetilde{V} = V \ast V $ and $ \widetilde{V}_{\epsilon} := V_{\epsilon} \ast
%V_{\epsilon } $

\begin{lemm} \label{PaleyWiener} Let $ g \in L^{\infty}(\Ra,\Ra) $ be supported in $ [0,b]
$, $  b > 0$, but in no smaller interval. Setting, for all $ 0 <
b^{\prime} < b $ and $ h \in (0,1] $,
$$ s_{b^{\prime}}(h) := \sup_{ 1 \leq |\xi| \leq 2 \atop \ {\rm Im}(\xi)<
0, \ {\rm Re}(\xi) \geq 0} | e^{ b^{\prime}{\rm Im}(\xi)/h}
( {\mathcal F}_{\rm inv} g) (\xi/h) |  ,
$$
we have
$$ \limsup_{h \rightarrow 0} s_{b^{\prime}}(h) = + \infty . $$
\end{lemm}

\noindent {\it Proof.} We clearly have
$$ | ( {\mathcal F}_{\rm
inv} g) (\xi) | \leq \frac{b}{2 \pi} ||g||_{\infty} e^{b |{\rm
Im}(\xi)|}, \qquad \xi \in \Ca, \ \mbox{Im}(\xi) \leq 0 , $$ and $
( {\mathcal F}_{\rm inv} g) $ is bounded for $
\mbox{Im}(\xi) > 0 $. Fix $ 0 < b^{\prime} < b $. By the
Paley-Wiener Theorem, we have
\begin{eqnarray}
 \sup_{ {\rm Im}(\xi)< 0 } | e^{- b^{\prime} |{\rm Im}(\xi)|}
( {\mathcal F}_{\rm inv} g) (\xi) | = + \infty ,
\label{supremum}
\end{eqnarray}
otherwise $ g $ should be supported in $ [0,b^{\prime}] $ which is
excluded. Furthermore, since $ g  $ is real valued, we have
$$  \left| ( {\mathcal F}_{\rm
inv} g) (\mbox{Re}(\xi)+i \mbox{Im}(\xi)) \right| =  \left|
( {\mathcal F}_{\rm inv} g) ( -\mbox{Re}(\xi)+i
\mbox{Im}(\xi)) \right|,
$$
so the supremum in (\ref{supremum}) can be taken over  $
\mbox{Re}(\xi) \geq 0 $ and $ \mbox{Im}(\xi)<0 $. Then, using the
local boundedness of $ ( {\mathcal F}_{\rm inv} g) $
and by writing the set $ \{ \xi \ | \ \mbox{Im}(\xi)< 0, \
\mbox{Re}(\xi) \geq 0 \} $ as
$$ \{ \xi \ | \ \mbox{Im}(\xi)< 0, \ \mbox{Re}(\xi) \geq 0 , \
|\xi|<1 \} \sqcup_{k \geq 0} \{ \xi \ | \ \mbox{Im}(\xi)< 0 , \
\mbox{Re}(\xi) \geq 0, \ 2^k \leq |\xi| < 2^{k+1} \} ,
$$
we have
$$ \limsup_{k \rightarrow + \infty} \sup_{  2^k \leq |\xi| < 2^{k+1}, \atop {\rm Im} (\xi)<0 , \ {\rm Re}(\xi) \geq 0 }
 | e^{- b^{\prime} |{\rm Im}(\xi)|} ( {\mathcal F}_{\rm
inv} g) (\xi) | = + \infty , $$ and the result follows.
 \finpreuve

\bigskip

\noindent {\bf Proof of Theorem \ref{theoremeprincipal3}.}
%Consider
%$$ W := \{  |z|e^{-i \theta} \ ; \ |z| \in [1,4] ,
%\ \theta \in [0,\pi] \} $$ and let $ \Omega $ be a simply
%connected neighborhood of $ W $ such that, for some $ \theta_0 <
%\pi $, (\ref{localiseOmega}) and (\ref{intervalle}) hold.
Fix  $ V \in C_0^{\infty}(\Ra,\Ra) $ with $ V \ne 0 $. By Theorem
\ref{theoremeprincipalp}, we can write
$$ D^{\zeta}_2 (H_0,H_0 + V ;
z,h ) = \prod_{w \in {\rm Res}(H_0+V, \Omega)} (z - w) \times \exp
(\varphi_2 (z,h))  $$ where, by Theorem \ref{theoremeprincipal2},
\begin{eqnarray}
 | \partial_z  \varphi_2 (z,h) | \lesssim h^{-1}, \qquad z \in W
. \label{borneabsorbable}
\end{eqnarray}
 On the other hand, by (\ref{coincidencedefinition}), $
D^{\zeta}_p (H_0,H_0 + V ; z,h ) $ can be replaced by the
definition (\ref{defFredholm}) using Fredholm determinants. Thus,
by (\ref{avecfacteurcorrectif}), we have
$$  D_3^{\zeta} (H_0,H_0 + V ; z,h) =  \prod_{w \in {\rm Res}(H_0+V, \Omega)} (z - w) \times  e^{\varphi_3(z,h)} $$
where
$$ \varphi_3 (z,h) := \varphi_2 (z,h) + \frac{\phi (z,h)}{2} , $$
with $ \phi $ given by (\ref{definitionphi}). In particular we
have
$$ ( \partial_z \phi )(z,h) = - \frac{ \pi}{2 z^{3/2} h^3}
(\partial_{\xi} {\mathcal F}_{\rm inv} \widetilde{V}^+_{\rm ev})
(2 z^{1/2} h^{-1}) + \frac{ \pi}{2 z^2 h^2}({\mathcal F}_{\rm inv} \widetilde{V}^+_{\rm ev})(2 z^{1/2} h^{-1})= \frac{ \pi}{2 z^2 h^2}f(2 z^{1/2} h^{-1}) , $$
where
$$f(\xi):= ({\mathcal F}_{\rm inv} \widetilde{V}^+_{\rm ev})(\xi) - \frac12 \xi \partial_{\xi}( {\mathcal F}_{\rm inv} \widetilde{V}^+_{\rm ev})(\xi)
=  ({\mathcal F}_{\rm inv} g)(\xi), $$
with 
$$g(x) := {\bf{1}}_{[0,+\infty)}\Big( \frac32   \left( \widetilde{V} (x) + \widetilde{V} (-x) \right)+  \frac12  x\partial_x \left( \widetilde{V} (x) + \widetilde{V} (-x) \right) \Big).$$
 Since $ V \ne 0 $, we have $
\widetilde{V}(0) = \int V^2
> 0 $ so $ g $ is supported in an
interval $ [0,b] $, $ b > 0 $, and no smaller one. We then obtain
(\ref{bornepresqueexponentielle}) with $ \delta = b/4 $,  first
remarking that, by (\ref{borneabsorbable}),
$$ | e^{\delta {\rm Im}(z^{1/2})/h} h\partial_z \varphi_2 (z,h) | \lesssim 1 , $$
 secondly that
\begin{eqnarray*}
 | e^{\delta {\rm Im}(z^{1/2})/h} h\partial_z \phi (z,h) | & \gtrsim & | h^{-1} e^{2 \delta {\rm Im(\xi)}
 / h^{\prime} }
 ( {\mathcal F}_{\rm
inv} g) (\xi / h^{\prime})
  |, \qquad \xi = z^{1/2}, \ \ h^{\prime} = h / 2 ,
\end{eqnarray*}
 and finally using Lemma \ref{PaleyWiener} with $ b^{\prime} = b/2 $. \finpreuve

\section{Analytic perturbations} \label{sectionanalytique}
In this section, we briefly prove a result similar to Theorem \ref{theoremeprincipal2} for $ p \geq 3 $
 in the more restrictive situation of analytic
perturbations. Namely, we consider $ V $ with coefficients analytic close to $ x = 0 $ (uniformly bounded with respect to $h$) and such that $ V \in {\mathcal
V}_{\rho}(\theta_0,R_0,\epsilon_0) , $ for any $R_0>0$. We denote by ${\mathcal
V}_{\rho}(\theta_0,0,\epsilon_0)$ the set of such perturbations $V$ and we assume that $ 0 < \theta_0 < \pi/2 . $
%Note that the elements of $ {\mathcal C}_{\rho}
%(\theta_0,0,\epsilon_0) $ are analytic close to $ x = 0 $ and thus
%so are the coefficients of $ V $. 
Here $ \rho > 0 $ is arbitrary.

In the following lemma, we first check that we can approximate such operators by fast
decaying ones. To avoid any confusion with $ \scal{x} = ( 1 + |x_1|^2 +
\cdots + |x_d|^2 )^{1/2} $, we set
$$ \langle \scal{x} \rangle =  (1 + x_1^2 + \cdots + x_d^2 )^{1/2},
\qquad \mbox{for} \   x \in \Ca^d \ \mbox{such that} \ 1 + x_1^2 +
\cdots + x_d^2 \notin (-\infty , 0 ] , $$ using the principal
determination of the square root mapping $ \Ca \setminus (-\infty
, 0 ] $ into $ e^{i (- \pi/2, \pi/2 )}(0,+\infty) $.
\begin{lemm} \label{approxdilate}  Let
$$ \chi_n (x)  = \exp \left( -   \langle \scal{x} \rangle / n \right),
  \qquad n \gg 1 , \ x \in \Ra^d .
$$ If $ \epsilon_0 $ is small enough, then, for $ n \geq n_0  $ large enough,
$$ V_n := \chi_n V \chi_n $$
belongs to $ {\mathcal V}_{\overline{d}} (\theta_0,0,\epsilon_0) $ for all $ \overline{d}
> d $,  the sequence $ (V_n)_{n \geq n_0} $ is bounded in $ {\mathcal V}_{\rho} (\theta_0,0,\epsilon_0) $
and, for all $ \rho^{\prime} < \rho $ and all $ s,\sigma \in \Ra
$,
\begin{eqnarray}
 || V_n - V ||_{H_{\rm sc}^{s,\sigma} \rightarrow H_{\rm
sc}^{s-2,\sigma+\rho^{\prime}}} \rightarrow 0, \qquad n
\rightarrow \infty , \label{compacitelimitedilate}
\end{eqnarray}
 for all $ h \ll 1 $.
\end{lemm}

\noindent {\it Proof.} The proof is similar to the one of Lemma
\ref{approxdistordu} (and anyway fairly elementary). The only new
point to check is that  the coefficients of $ V_n $ belong to $
{\mathcal V}_{ \overline{d} } (\theta_0,0,\epsilon_0) $ and are
bounded in $ {\mathcal V}_{\rho} (\theta_0,0,\epsilon_0) $.
Indeed,  for $ r = e^{i \theta} t $, with $ t > 0 $ and $ \theta
\in [0,\theta_0] $, and for $ \omega $ such that $
\mbox{dist}_{\Ca^d}(\omega,{\mathbb S}^{d-1}) < \epsilon_0 $, we
first note that, if $ \epsilon_0 $ is small enough, $ r^2 \omega^2
\notin (-\infty,0 ] $. Furthermore, if $t$ is large, $ 1 + r^2
\omega^2 = t^{2} e^{2 i \theta} \left(1 + o(1) \right) $, thus
$$  {\rm Re} \langle \scal{r \omega} \rangle  \gtrsim t \cos (\theta) . $$
It is then easy to check that, for all $ \alpha  $, $
\partial^{\alpha} \chi_n $ is bounded on $  \Sigma (\theta_0,0, \epsilon_0 ) $, uniformly with respect
to $n \geq 1$.
%$$ \chi_n (r \omega) = \exp ( -  r^2 \omega^2 / n ) =
%\exp \left( - \cos (2 \theta ) t^2 \omega^2/n  \right) \exp \left(
%- i \sin (2 \theta ) t^2 \omega^2/n \right) .
%$$
%Since $ | 2\theta  | < 2 \theta_0 < \pi /2  $, there exists $
%\gamma
%> 0 $ such that $ \cos (2 \theta ) \geq 2 \gamma $ and
%thus, if $ \epsilon_0 $, is small enough, we have
%$$ | \exp ( -  r^2 \omega^2 / n ) | \leq C
%\exp \left( - \gamma t^2  /n  \right) . $$ This exponential decay
% then shows that, if $ \epsilon_0 $ is small enough,
%\begin{eqnarray}
% \left| \exp ( - 2 \cos ( 2 \theta  ) x^2 / n ) x^{\alpha}/
%n^{|\alpha|} \right| \leq C_{\alpha} \exp \left( - \gamma t^2  /n
%\right) , \label{decroissanceexponentielleexplicite}
%\end{eqnarray}
% for all $n \geq 1 $ and all $ x = e^{i \theta } t
%\omega $ with $$  n \geq 1, \ t > 0 , \ \theta \in [0,\theta_0], \
%\mbox{dist}_{\Ca^d}(\omega,{\mathbb S}^{d-1}) < \epsilon_0 . $$
%Note the uniformity with respect to $ n \geq 1 $.
 Since the
coefficients of $ V_n $ are linear combinations of products of
coefficients of $ V $ by $ \chi_n
\partial_{x}^{\alpha}  \chi_n   $,
%which behave like (\ref{decroissanceexponentielleexplicite}),
we see that $ (V_n)_{n \geq 1}  $ is  bounded in $ {\mathcal
V}_{\rho} (\theta_0,0,\epsilon_0) $. It also clearly belongs to
$\in {\mathcal V }_{\overline{d}} (\theta_0,0,\epsilon_0) $.
\finpreuve

\bigskip

We next give an elementary deformation result along $ e^{i \theta
 }\Ra^d $. Let us denote
 $$ V_{\rm dil}(\theta) := \sum_{|\alpha| \leq 2} v_{\alpha} (e^{i \theta}x,h) (e^{- i  \theta}h D)^{\alpha} , $$
if $ V = \sum_{|\alpha| \leq 2} v_{\alpha} (x,h) (h D)^{\alpha} $, that
is (\ref{parametrisation}) with $ \kappa (x) = e^{i \theta} x $ and $P=V$.
For $ i \theta \in \Ra $, we also have
$$ V_{\rm dil}(\theta) = U_{\rm dil}(i \theta) V U_{\rm dil}(i \theta)^* ,  $$
where $ U_{\rm dil}(t) $ is the generator of dilations introduced
for similar purposes in \cite{AguilarCombes}
$$ U_{\rm dil}(t) u (x) = e^{td/2} u (e^{t} x) . $$
\begin{lemm} \label{lemmededilatation} Let $ k > d/2 + 1 $. For all $ n \gg 1 $, $ \theta \in [0,\theta_0] $,
$ z \in \Omega^+ $ and $ j \geq 1 $,
\begin{eqnarray}
 \emph{tr} \left( \frac{d^j}{d \varepsilon^j} (H_0 + \varepsilon
V_n - z)^{-k}_{| \varepsilon = 0} \right) = \emph{tr} \left(
\frac{d^j}{d \varepsilon^j} (e^{-2 i \theta }H_0 + \varepsilon
V_{n, {\rm dil}}(\theta) - z)^{-k}_{| \varepsilon = 0} \right) .
\label{dilatationtracederivee}
\end{eqnarray}
\end{lemm}

\noindent {\it Proof.} For $ i \theta \in \Ra $, the result is
obvious since the right hand side of
(\ref{dilatationtracederivee}) reads
$$ \mbox{tr} \left( \frac{d^j}{d \varepsilon^j} U_{\rm dil}(i \theta) (H_0
+ \varepsilon V_{n} - z)^{-k} U_{\rm dil}(i \theta)^*_{|
\varepsilon = 0} \right) .
$$
On the other hand, $ \theta \mapsto V_{n,{\rm dil}}(\theta) $ is
holomorphic from $ (0,\theta_0) + i (-1,0) $ to $ {\mathcal
L}(H_{\rm sc}^{s+2 , \sigma},H_{\rm sc}^{s,\sigma + \overline{d}})
$, for all $s \in \Na, \sigma \in \Ra $ and $ \overline{d}>d $. It
is also continuous for $ \theta \in [0,\theta_0] + i [-1,0] $.
Since $ e^{-2i\theta}H_0 - z $ is invertible,  Proposition
\ref{lemmeprop} proves the existence of the resolvent $ (e^{-i2
\theta }H_0 + \varepsilon V_{n,{\rm dil}}(\theta) - z)^{-1} $ for $
\varepsilon  $ small enough (depending on $h$ but this harmless for we shall eventually set $ \varepsilon = 0 $). It is then
holomorphic for $ \theta \in (0,\theta_0) + i (-1,0) $ and
continuous for $ \theta \in [0,\theta_0] + i [-1,0] $, with values
in $ {\mathcal L}(H_{\rm sc}^{s , \sigma},H_{\rm sc}^{s+2,\sigma})
$. Therefore the expression of the right hand side of
(\ref{dilatationtracederivee}) given by Proposition
\ref{pourlimite} is holomorphic with respect to $ \theta \in
(0,\theta_0) + i (-1,0) $, continuous on $ [0,\theta_0] + i [-1,0]
$ and constant on $ i [-1,0] $ hence constant in $ [0,\theta_0] +
i [-1,0] $ by analytic continuation. This completes the proof.
\finpreuve
% in a connected neighborhood
%of $ [0,\theta_0] $. To check this holomorphy, we use that
%$$ \partial_{\theta} v_n (e^{i \theta}x,h) = i e^{i \theta}x \cdot \partial_x v_n (e^{i \theta} x,h)  $$
%decays fast (exponentially) at infinity and that
%$$ \partial_{\theta} (e^{-2 i \theta}H_0 -z)^{-1} = -2 i e^{-2 i \theta_0} H_0  (e^{-2 i \theta}H_0 -z)^{-2} . $$
% The result then follows by analytic continuation.

\bigskip

Next, using Propositions \ref{inverse}, \ref{lemmeprop}, \ref{pourlimite}, Lemma \ref{approxdilate} and the notation (\ref{Taylorkp}), we can write, for each $ z \in \Omega^+
$,
\begin{eqnarray}
 \zeta_p (k,z,h) & = & \lim_{n \rightarrow \infty} T_p^k
(H_0,V_n,z) , \nonumber
\end{eqnarray}
that is the limit of
\begin{eqnarray}
\mbox{tr} \left( (H_0 + V_n-z)^{-k} - (H_0-z)^{-k} \right) -
\sum_{j=1}^{p-1} \frac{1}{j !} \mbox{tr} \left( \frac{d^j}{d
\varepsilon^j} (H_0 + \varepsilon V_{n} - z)^{-k}_{| \varepsilon =
0} \right) , \nonumber
\end{eqnarray}
or, by Lemma \ref{lemmededilatation}, the limit of
\begin{eqnarray}
\mbox{tr} \left( (H_0 + V_n-z)^{-k} - (H_0-z)^{-k} \right) -
\sum_{j=1}^{p-1} \frac{1}{j !} \mbox{tr} \left( \frac{d^j}{d
\varepsilon^j} (e^{-2 i \theta_0 }H_0 + \varepsilon V_{n, {\rm
dil}}(\theta_0) - z)^{-k}_{| \varepsilon = 0} \right) . \nonumber
\end{eqnarray}
Observing that Proposition \ref{sjos2} can be extended to the
sequence $ V_n $ (ie that the corresponding finite rank operators
$ K_{n}(\theta_0) $ converge as $ n \rightarrow + \infty $), this
limit is the sum of
\begin{eqnarray}
  \mbox{tr} \left( (\widehat{H_1(\theta_0)} - z)^{-k} -
(\widehat{H_0(\theta_0)} - z)^{-k} - \sum_{j=1}^{p-1} \frac{1}{j
!} \frac{d^j}{d \varepsilon^j} (e^{-2 i \theta_0 }H_0 +
\varepsilon V_{ {\rm dil}}(\theta_0) - z)^{-k}_{| \varepsilon = 0}
\right) , \label{operateureneffetatrace}
\end{eqnarray}
 and of
\begin{eqnarray}
- \frac{\partial_z^{k-1}}{(k-1)!} \left(
 \partial_z \log \mbox{det}_1 \left( 1 +
  \widetilde{K}_{1} (\theta_0,z) \right) - \partial_z \log \mbox{det}_1 \left( 1 +
  \widetilde{K}_0 (\theta_0,z) \right) \right)  = \nonumber \\
 \qquad  \qquad \sum_{w \in {\rm Res}(H_0+V,\Omega)}
\frac{1}{(w-z)^k} + \phi(z,h), \nonumber
%\label{contributionresonances}
\end{eqnarray}
with $ \phi(z,h) $  holomorphic on $ \Omega $ and $ {\mathcal
O}(h^{-d}) $ locally uniformly. This follows from
 (\ref{ajoutcompact}), (\ref{deriveelogarithmique}),
(\ref{factorisationSjostrand}), Proposition \ref{Harnack} and from
the absence of resonances for $ H_0 $. The operator inside the
trace in (\ref{operateureneffetatrace}) is  trace class because
it is the sum of
 of
\begin{eqnarray}
 (\widehat{H_1(\theta_0)} - z)^{-k} - \sum_{j=0}^{p-1} \frac{1}{j
!} \frac{d^j}{d \varepsilon^j} (e^{-2 i \theta_0 }H_0 +
\varepsilon V_{ {\rm dil}}(\theta_0) - z)^{-k}_{| \varepsilon =
0} , \label{modulocompact}
\end{eqnarray}
and of
$$ (\widehat{H_0(\theta_0)} - z)^{-k}  -  (e^{-2i \theta_0}H_0 - z)^{-k} $$
which is $ {\mathcal O}(h^{-d}) $ in the  trace class for $z \in
\Omega$ by Propositions \ref{inverse}, \ref{pourlimite} (recall that $ \widehat{H_0(\theta_0)}  -  e^{-2i \theta_0}H_0 $ is compactly supported) and \ref{sjos2}, using the elementary bound $ || (e^{-2i\theta_0}H_0 - z)^{-1}||_{L^2 \rightarrow H^{2,0}_{\rm sc}} \lesssim 1 $.
Setting
$$ \widehat{V}(\theta_0) := \widehat{H_1(\theta_0)}  - (  e^{-2 i \theta_0 }H_0 +  V_{
{\rm dil}}(\theta_0) ) $$ which is compactly supported, 
(\ref{modulocompact}) is the sum of the trace class operators 
$$ \frac{1}{j!} \frac{d^j}{d \varepsilon^j} \left( (e^{-2 i \theta_0 }H_0 +
\varepsilon V_{ {\rm dil}}(\theta_0) +
\varepsilon  \widehat{V}(\theta_0)  - z
)^{-k}  -   (e^{-2 i \theta_0 }H_0 +
\varepsilon V_{ {\rm dil}}(\theta_0) - z)^{-k} \right)_{| \varepsilon =
0}  $$
and of
%which are holomorphic with respect to $ z \in \Omega $ and with trace norm locally uniformly of order $ h^{-d} $, and of
$$  ( e^{-2 i \theta_0 }H_0 +
 V_{ {\rm dil}}(\theta_0) +
 \widehat{V}(\theta_0)  - z)^{-k} - \sum_{j=0}^{p-1} \frac{1}{j
!} \frac{d^j}{d \varepsilon^j} \left( e^{-2 i \theta_0 }H_0 +
\varepsilon V_{ {\rm dil}}(\theta_0) +
\varepsilon  \widehat{V}(\theta_0)  - z
\right)^{-k}_{| \varepsilon = 0} , $$ 
which are all of order $ h^{-d} $ in the  trace
class, locally uniformly with respect to $z \in \Omega$, by Proposition \ref{pourlimite}, (\ref{bornefinale}), Proposition \ref{sjos2}
and again the estimate $ || (e^{-2i\theta_0}H_0 - z)^{-1}||_{L^2 \rightarrow H^{2,0}_{\rm sc}}  \lesssim 1 $.

Using Proposition \ref{reductionZetaintro}, we obtain the
following theorem.

\begin{theo} Let $ \rho > 0 $ and $ p \in \Na $ such that $ p \rho > d $.
Let $ \Omega \Subset e^{- i(2 \theta_0,\epsilon)} (0, + \infty) $
be a simply connected open subset with $ 0 < \theta_0 < \pi/2 $, $
\epsilon > 0 $ small enough and satisfying (\ref{intervalle}).
Then, if $ V \in {\mathcal V}_{\rho}(\theta_0,0,\epsilon_0) $,
% for all $ R_0
%> 0 $ with coefficients analytic close to $ 0 $ and 
%bounded   uniformly with respect to $h$,  
any
$ \varphi_p $ as in Theorem \ref{theoremeprincipalp} satisfies,
for all $ W \Subset \Omega $,
$$ | \partial_z \varphi_p (z,h) | \leq C_W h^{-d}, \qquad z \in W , \ h \ll 1 . $$
\end{theo}

\bigskip

\noindent {\bf Acknowledgement.} We are pleased to dedicate this paper to Didier Robert. It was started on the occasion of his 60th anniversary and answers a question he raised a few years ago.


\begin{thebibliography}{99}

\bibitem{AguilarCombes} {\sc J. Aguilar, J.-M. Combes}, {\it A class of analytic perturbations for
one-body Schr\"odinger Hamiltonians}, Comm. Math. Phys. { 22},
269--279  (1971).

\bibitem{Belov} {\sc S. M. Belov, A. V. Rybkin}, {\it Higher order trace formulas of the Buslaev-Faddeev type for
the half line Schr\"odinger operator with long range potentials}, J. Math. Phys. 44, 7, 2748-2761 (2003).

\bibitem{BiKr} {\sc M. Sh. Birman, M. G. Krein}, 
{\it On the theory of wave operators and scattering operators}, Sov. Math. Dokl. 3, 740-744 (1962).

\bibitem{BoBrRa} {\sc J.-F. Bony, V. Bruneau and  G. Raikov},
{\it  Resonances and spectral shift function near the landau
levels},  Ann. Inst. Fourier { 57},  no. 2, 629-671  (2007).


\bibitem{BoJuPe} {\sc D. Borthwick, C. Judge, P.A. Perry},
{\it Determinants of Laplacians and isopolar metrics on surfaces
of infinite area},  Duke Math. J. { 118},    no. 1, 61-102 (2003).


\bibitem{BoucletAA} {\sc J.-M. Bouclet},
{\it Trace formulae for relatively Hilbert-Schmidt perturbations}, Asymp. Analysis 32, no. 3-4, 257-291 (2002).



\bibitem{BoucletJFA} {\sc \name}, {\it Spectral distributions for long range perturbations}, J.
Funct. Anal. 212, no. 2, 431-471 (2004).


\bibitem{BruneauPetkov} {\sc V. Bruneau, V. Petkov}, {\it Meromorphic continuation of the spectral shift
function}, Duke Math. J. vol. 116, No. 3, 389-430 (2003).

\bibitem{Carr1} {\sc G. Carron},
{\it D\'eterminant relatif et la fonction Xi,} Amer. J. Math, vol.
124  n. 2, 307-352 (2002).

\bibitem{DiSj} {\sc M. Dimassi, J. Sj\"ostrand}, {\it Spectral asymptotics in the
semi-classical limit}, London Mathematical Society, Lecture Notes
Series { 268}, Cambridge University Press (1999).


\bibitem{Froese1} {\sc R. Froese},  {\it Asymptotic distribution of resonances in one dimension},
J. Diff. Equa. { 137},  251-272 (1997).


\bibitem{Froese2} {\sc \name}, {\it Upper bounds for the resonance counting
function of Schr\"odinger operators in odd dimensions}  Canad. J. Math. { 50}
  no. 3, 538-546 (1998). Correction in  Canad. J. Math. { 53} ,  no. 4, 756-757 (2001).

\bibitem{GPS} {\sc F. Gesztesy, A. Pushnitski, B. Simon}, {\it On the Koplienko spectral shift function, I. Basics},
Journal of Mathematical Physics, Analysis, Geometry, to appear. 


\bibitem{GoKr1}{\sc I. C. Gohberg, M. G. Krein},
{\it Introduction to the theory of linear nonselfadjoint
operators,} Translations of Mathematical Monographs, Vol. 18. AMS.
Providence, R.I. (1969).


\bibitem{Gu06} {\sc C. Guillarmou}, {\it Generalized Krein formula and determinants for even dimensional Poincar\'e-Einstein manifolds}, preprint, arXiv:math/0512173.


\bibitem {GuZw} {\sc L. Guillop\'e, M. Zworski}, {\it Scattering asymptotics for Riemann surfaces},
Ann. of Math. (2) 145, n. 3, 597-660 (1997).

\bibitem{Krei1} {\sc M. G. Krein},
{\it Perturbation determinants and a formula for the trace of
unitary and selfadjoint operators}, Soviet Math. Dokl. 3, 707-710
(1962).


\bibitem{Horm3} {\sc L. H\"ormander}, {\it The Analysis of Linear PDO III},  Springer Verlag (1985).

\bibitem{Koplienko}{\sc L. S. Koplienko}, {\it Trace formula for non trace class
perturbations}, Siberian Math. J. 25, 735-743 (1984).

\bibitem{Koplienko2} {\sc \name}, {\it Regularized spectral shift function for one dimensional Schr\"odinger operator with slowly decreasing
potential}, Siberian Math.  J. 26, 365-369 (1985).

\bibitem{Mart} {\sc A. Martinez}, {\it An introduction to semiclassical and microlocal
analysis}, Universitext, Springer-Verlag, New-York (2002).

\bibitem{Muller} {\sc W. M\"uller}, {\it Relative zeta functions, relative
determinants and scattering theory}  Comm. Math. Phys.  { 192},
no. 2, 309-347  (1998).

\bibitem{PeZw99} {\sc V. Petkov, M. Zworski}, {\it Breit-Wigner approximation and the distribution of resonances},
  Comm. Math. Phys. {  204},  no. 2, 329-351 (1999). Erratum:  Comm. Math. Phys.   { 214},  no. 3, 733-735  (2000).

\bibitem{PeZw01} {\sc \name}, {\it Semi-classical estimates on the scattering determinant},
 Ann. Henri Poincar\'e { 2},  no. 4, 675-711   (2001).


\bibitem{RaySinger} {\sc D. B. Ray, I. M. Singer}, {\it R-Torsion and the Laplacian on Riemannian
manifolds}, Adv. in Math. 7, 145-210 (1971).



\bibitem{RoBook} {\sc D. Robert},  {\it Autour de l'approximation
semi-classique}, Progress in mathematics, {68}, Birkha\"user, Boston
(1987).

\bibitem{Rybkin} {\sc A. V. Rybkin}, {\it On a trace formula of the Buslaev-Faddeev type for a long range potential}, J. Math. Phys 40, 3, 1334-1343 (1999).


\bibitem{Simon} {\sc B. Simon}, {\it Resonances in one dimension and Fredholm
determinants}, J. Funct. Anal. 178, 396-420 (2000).

\bibitem{Sjostrace} {\sc J. Sj\"ostrand}, {\it A trace formula and review of some estimates for
resonances}, p. 377-437 in Microlocal Analysis and Spectral
Theory, NATO ASI Series C, vol. 490, Kluwer (1997).

\bibitem{Sjosequi} {\sc \name}, {\it Resonances for bottles and trace formulae}, Math. Nach. { 221}, 95-149 (2001).

\bibitem{SjBook} {\sc \name}, {\it Lectures on resonances},
available at: www.math.polytechnique. fr/$\sim$sjoestrand/.

\bibitem{SjostrandZworski} {\sc J. Sj\"ostrand, M. Zworski}, {\it Complex scaling and the distribution of scattering
poles}, J. Amer. Math. Soc. 4, n. 4 729-769 (1991).

\bibitem{Vodev} {\sc G. Vodev}, {\it Sharp bounds on the number of scattering poles for perturbations of the Laplacian},
Comm. Math. Phys. 146, 205-216 (1992).

\bibitem{Watson} {\sc G.N. Watson}, {\it A treatise on the theory of Bessel functions}, Cambridge University Press, Cambridge (1966).

\bibitem{Yafaev} {\sc D. Yafaev}, {\it Mathematical Scattering Theory}, AMS,
Providence, RI (1992).

\bibitem{Zw89} {\sc M. Zworski}, {\it  Sharp polynomial bounds on the number of scattering poles},
  Duke Math. J.  { 59},  no. 2, 311-323 (1989).


\end{thebibliography}
\end{document}